\documentclass[10pt,twoside]{amsart}
\usepackage{amsmath, amsfonts, amsthm, amssymb, verbatim}
\usepackage[mathcal]{euscript}
\newcommand \RR   {\mathbb{R}}

\newcommand \be  {\begin{equation}}
\newcommand \ee  {\end{equation}}
\newcommand \bs  {\begin{split}}
\newcommand \es  {\end{split}}
\newcommand \del {{\partial}}
\newcommand \eps {\epsilon}

\newcommand \la{\langle}
\newcommand \ra{\rangle}
\newcommand \supp {\mbox{supp }}
\newcommand \sgn {\mbox{sgn}}
\newcommand \overf {{\overline f}}

\newcommand \DD    {\mathcal{D}}
\newcommand \1    {\mathbf{1}}
\newcommand \LLL    {\mathbf{L}}
\newtheorem{definition}{Definition}[section]
\newtheorem{lemma}[definition]{Lemma}
\newtheorem{theorem}[definition]{Theorem}

\numberwithin{equation}{section}

\begin{document}
\bibliographystyle{plain}
\title[Isothermal Compressible Fluids]
{Symmetries and global solvability 
\\
of the isothermal gas dynamics equations
}
\author[LeFloch and Shelukhin]
{\scshape
Philippe G. LeFloch$^1$ \and Vladimir Shelukhin$^2$}
\thanks
{$^1$ Laboratoire Jacques-Louis Lions 
\& Centre National de la Recherche Scientifique, University of Paris 6,
4 place Jussieu,  75252 Paris, France. 
E-mail : {\sl lefloch@ann.jussieu.fr}} 
\thanks{
$^2$ Lavrentyev Institute of Hydrodynamics,
Prospect Lavrentyeva 15, Novosibirsk, 630090, Russia.
E-mail : {\sl shelukhin@hydro.nsc.ru}
\newline
2000\textit
{\ AMS Subject Classification:} 35L, 35L65, 76N, 76L05
\newline
\textit{Key Words:} Euler equations, isothermal compressible fluids,
mathematical entropy, compensated compactness, existence theory. 
}
\begin{abstract}
We study the Cauchy problem associated with the system of two conservation laws
arising in isothermal gas dynamics, in which
the pressure and the density are related by the $\gamma$-law equation
$p(\rho) \sim \rho^\gamma$ with $\gamma =1$.
Our results complete those obtained earlier for $\gamma >1$.
We prove the global existence and compactness of entropy solutions generated by
the vanishing viscosity method. The proof relies on compensated
compactness arguments and symmetry group analysis. Interestingly,
we make use here of the fact that the isothermal gas dynamics
system is invariant modulo a linear scaling of the density.
This property enables us to reduce our problem to that with a small initial density.

One symmetry group associated with the linear hyperbolic equations
describing all entropies of the Euler equations gives rise to a
fundamental solution with initial data imposed to the line
$\rho=1$. This is in contrast to the common approach (when $\gamma
>1$) which prescribes initial data on the vacuum line $\rho =0$. The
entropies we construct here are weak entropies, i.e. they vanish
when the density vanishes. Another feature of our proof lies in the reduction theorem which
makes use of the family of weak entropies to show that a Young measure
must reduce to a Dirac mass. This step is based on new convergence results
for regularized products of measures and functions of bounded variation.
\end{abstract}
\maketitle



\section{Introduction}
\setcounter{equation}{0}

We consider the Euler equations for  compressible fluids
\begin{align}
&  \del_t \rho + \del_x (\rho u) = 0,
\label{1.1}
\\
&  \del_t (\rho u) + \del_x( \rho u^2 + p(\rho)) = 0,
\label{1.2}
\end{align}
where $\rho\geq 0$ denotes the density, $u$ the velocity, and
$p(\rho) \geq 0$ the pressure. We assume that the fluid
is governed by the isothermal equation of state
\be
p(\rho)
= k^2 \, \rho,
\label{1.3}
\ee
where $k>0$ is a constant. Observe that the scaling $u \to k \, u$, $t\to t/k$
allows one to reduce the system (1.1)--(1.3) to the same system with $k=1$.

The existence of weak solutions (containing jump discontinuities) for the Cauchy problem
associated with (1.1)--(1.3) was first established by Nishida \cite{Nishida}
(in the Lagrangian formulation). The solutions obtained by Nishida
have bounded variation and remain bounded away from the vacuum. For background
on the BV theory we refer to \cite{Dafermos,LeFloch2}.

By contrast, we are interested here in solutions in a much weaker functional class
and in solutions possibly reaching the vacuum $\rho=0$. Near the vacuum, the system
(1.1)--(1.3) is degenerate and, in particular, the velocity $u$ can not be defined uniquely.
Indeed, the present paper is devoted to developing the existence
theory in
a framework covering solutions satisfying
$$
\rho \in L^\infty (\Pi), \quad
\rho |u| \leq C \, ( \rho + \rho |\log \rho|),\quad
\Pi =\RR \times (0,T),
$$
with a constant $C>0$ depending solely on initial data.
The time interval $(0,T)$
is arbitrary. Our proof extends DiPerna's pioneering work \cite{DiPerna2}
concerned with the pressure law $p(\rho) \sim \rho^\gamma$.


\section{Main result}
\setcounter{equation}{0}

Introducing the momentum variable $m:=\rho u$, one can reformulate the Cauchy problem
associated with (1.1)--(1.3) as follows:
\be
\begin{split}
& \del_t \rho + \del_x m =0,
\\
& \del_t m + \del_x (\frac{m^2}{\rho} + \rho) = 0,
\end{split}
\label{2.1a}
\ee
with initial condition
\be
\begin{split}
& \rho|_{t=0} = \rho_0, \qquad m|_{t=0} = m_0 := \rho_0 \, u_0.
\end{split}
\label{2.1b}
\ee
where $\rho_0, u_0$ are prescribed. Let us first recall the following terminology.
A pair of (smooth) functions $\eta=\eta (m,\rho)$, $q=q(m,\rho)$ is called an {\sl entropy pair}
if, for any smooth solution $(m,\rho)$ of \eqref{2.1a}, one also has
$$
\del_t \eta (m,\rho)+\del_x q(m,\rho) = 0.
$$
More precisely, we consider entropies $\eta,q \in C^2 (\Omega) \cap C^1(\bar{\Omega})$
in any domain of the form
$$
\Omega := \bigl\{0<\rho <\rho_{\ast},\quad |m|<c_{\ast} \rho \, (1+|\ln\rho|) \bigr\},
\quad c_{\ast}>0,\quad\rho_{\ast}>0.
$$
It is easily checked that $\eta, q$ must solve the equations 
\be
q_m = 2 \, \frac{m}{\rho}\eta_m +\eta_{\rho},\quad
q_{\rho}=\eta_m -\frac{m^2}{\rho^2}\eta_m,
\label{2.2a}
\ee
which implies that
\be
\eta_{\rho\rho}=\frac{p' (\rho)}{\rho^2}\eta_{uu}
= \frac{1}{\rho^2}\eta_{uu}.
\label{2.2b}
\ee
A pair $(\eta, q)$ is said to be a {\sl weak entropy} if
$\eta(0,0)=q(0,0)=0$. It is said to be {\sl convex}
if in addition, $\eta$ is convex
with respect to the conservative variables $(\rho,m)$.

Given an initial data $m_0$, $\rho_0 \in L^\infty(\RR)$ obeying the inequalities
\be
\rho_0 (x)\geq 0,\quad
|m_0 (x)|\leq
c_0 \, \rho_0 (x) \, (1+|\ln\rho_0 (x)|), \qquad x \in \RR
\label{2.3}
\ee
for some constant $c_0 >0$, {\sl an entropy solution} to the Cauchy
problem \eqref{2.1a}-\eqref{2.1b} on the time interval $(0,T)$ is, by definition,
a pair of functions $(m,\rho)\in L^\infty(\Pi)$ satisfying
the inequalities
\be
\rho(x,t) \geq 0,\quad |m (x,t)|\leq c\rho (x,t)
(1+|\ln\rho  (x,t)|), \qquad (x,t) \in \Pi
\label{2.4}
\ee
for some positive constant $c$, together with the inequality
\be
\iint_\Pi  \Big( \eta(m,\rho) \, \del_t \varphi + q(m,\rho) \, \del_x \varphi \Big) \, dxdt
+ \int_\RR \eta(m_0, \rho_0) \, \varphi(\cdot,0) \, dx
\geq 0
\label{2.5}
\ee
for every convex, weak entropy pair $(\eta,q)$
and every non-negative function ${\varphi \in \DD (\RR\times [0,T))}$
(smooth functions with compact support).

\

The main results established in the present paper are summarized in Theorems 2.1--2.3 below.

\begin{theorem} {\rm (Cauchy problem in
momentum-density variables.)}
Given an arbitrary time interval $(0,T)$ and an initial data
$(m_0, \rho_0) \in L^\infty(\RR)$ satisfying the condition \eqref{2.3},
there exists an entropy solution $(m,\rho)$ of the Cauchy problem
\eqref{2.1a}-\eqref{2.1b} satisfying the inequalities \eqref{2.4},
with a constant $c$ depending on $c_0$ only.
\end{theorem}

To prove this theorem it will be convenient to introduce the Riemann invariants $W$ and $Z$
by
$$
W := \rho e^{u}, \quad Z := \rho e^{-u},
$$
or equivalently
$$
\rho =f_1 (W,Z):=(WZ)^{1/2},\qquad
\rho u =f_2 (W,Z):=(WZ)^{1/2}\ln (W /Z)^{1/2}.
$$
One can then reformulate the Cauchy problem \eqref{2.1a}-\eqref{2.1b} in terms of $W,Z$,
as follows
\be
\begin{split}
& \del_t f_1(W,Z) +\del_x f_2(W,Z) = 0,
\\
& \del_t f_2(W,Z) +\del_x (f_3(W,Z) + f_1(W,Z)) = 0,
\qquad
f_3 :=(WZ)^{1/2}\bigl(\ln (W /Z)^{1/2}\bigr)^2,
\end{split}
\label{2.6}
\ee
\be
W|_{t=0}=W_0 :=\rho_0 e^{u_0},\qquad
Z|_{t=0}=Z_0 :=\rho_0 e^{-u_0}.
\label{2.7}
\ee
a pair of non-negative functions $W,Z \in L^\infty(\Pi)$ if then called
an {\sl entropy solution} to the problem \eqref{2.6}-\eqref{2.7}
if
$$
\iint_\Pi \Big( \widetilde{\eta}(W,Z)\Big)\, \del_t \varphi
      + \widetilde{q}(W,Z) \, \del_x \varphi \Big) \, dxdt
+ \int_\RR \widetilde\eta(W_0, Z_0) \, \varphi(\cdot,0) \, dx
 \geq 0
$$
for any non-negative function $\varphi \in \DD (\RR\times [0,T))$,
where
$$
\tilde{\eta}(W,Z):=\eta (f_2 (W,Z),f_1 (W,Z)),\quad \tilde{q}(W,Z):=q(f_2 (W,Z),f_1 (W,Z)),
$$
and $(\eta, q)$ is any convex, weak entropy pair in the sense introduced above.

Theorem 2.1 above will be obtained as a corollary of the following result.

\begin{theorem}  {\rm (Cauchy problem in Riemann invariant variables.)}
Given non-negative functions $W_0,Z_0 \in L^\infty(\RR)$, the Cauchy problem
\eqref{2.6}-\eqref{2.7} has an entropy solution on any time interval (0,T).
\end{theorem}

It is checked immediately that, if $(W,Z)$ is an entropy solution
given by Theorem 2.2, then the functions $m:=f_2 (W,Z)$ and $\rho :=f_1 (W,Z)$
determine an entropy solution of the problem \eqref{2.1a}-\eqref{2.1b}.

One more consequence of Theorem 2.2 concerns the original problem \eqref{1.1}--\eqref{1.3}
in the density-velocity variables. Defining the density and velocity from the Riemann variables
by
$$
u := \ln (W/Z)^{1/2},
\qquad
\rho := (WZ)^{1/2},
$$
we deduce also the following result from Theorem 2.2.

\begin{theorem} {\rm (Cauchy problem in velocity-density variables.)}
Let $(0,T)$ be a time interval. Given any measurable functions $u_0$ and $\rho_0$
satisfying the conditions
$$
0 \leq \rho_0 \in L^\infty(\RR),
\qquad |u_0 (x)|\leq c_0 (1+|\ln \rho_0 (x)|), \quad x \in \RR
$$
for some positive constant $c_0$, there exist measurable functions
$u=u(x,t)$ and $\rho=\rho (x,t)$ such that
$$
0 \leq \rho \in L^\infty(\Pi),
\qquad |u(x,t)|\leq c (1+|\ln \rho (x,t)|), \quad (x,t) \in \Pi
$$
(where $c>0$ is a constant depending on $c_0$)
and $(u,\rho)$ is an entropy  solution of the problem \eqref{1.1}-\eqref{1.3}
in the sense that the entropy inequality
$$
\iint_\Pi  \Big(\eta(\rho,\rho \, u) \, \del_t \varphi
  + q(\rho, \rho \, u) \, \del_x \varphi \Big) \, dxdt
+ \int_\RR \eta(\rho_0, \rho_0 \, u_0) \, \varphi(\cdot,0) \, dx
\geq 0
$$
holds for any convex, weak entropy pair $(\eta,q)$
and any function $\varphi$ as in Theorem~2.1.
\end{theorem}

\

The novel features of our proof of the above results are :
\begin{itemize}
\item the use of symmetry and scaling properties
of both the isothermal
Euler equations and
the entropy-wave equation,
\item an analysis of new nonconservative products of functions with bounded variation by measures.
\end{itemize}
We rely on two classical ingredients. The first tool is the compensated compactness method
introduced by Tartar in \cite{Tartar1,Tartar2}. (See also Murat \cite{Murat1}.)
This method allows to show that a weakly convergent sequence (of approximate solutions given by the viscosity
method) is actually strongly convergent: such a result is achieved by a ``reduction lemma''
(to point mass measures) for Young measures representing the limiting behavior of the sequence.
Tartar method was applied to systems of conservation laws by DiPerna \cite{DiPerna1,DiPerna2}.
For a completely different approach to the vanishing viscosity method, we refer to 
Bianchini and Bressan \cite{BB}. Still another perspective is introduced in LeFloch \cite{LeFloch3}.

The second main tool is the symmetry group analysis of differential equations which goes back to
Lie's classical works. The first symmetry property we use concerns the system (1.1)--(1.3)
itself: we observe that it is invariant with respect to the scaling $\rho \to \lambda \, \rho$
($\lambda$ being an arbitrary parameter). This property allows us to assume that the density
is sufficiently small when performing the reduction of the Young measures.

To generate the class of weak entropies, we calculate all the Lie groups associated
with the entropy equation \eqref{2.2b} for the function $\eta$. By using one of them
we construct the fundamental solution with initial data prescribed on the line $\rho =1$.
This is in contrast with the standard approach which prescribes initial data on the vacuum line $\rho =0$.

The need of a large family of weak entropies for the Young measure reduction was demonstrated by DiPerna
for the isentropic gas dynamics equations with the pressure law $p=\rho^{\gamma}$, $\gamma >1$.
When $\gamma =\frac{2n+3}{2n+1}$, with $n$ being integer, DiPerna used weak entropies which are
progressive waves given by Lax. The method of Tartar and DiPerna was then 
extended by Serre \cite{Serre} to strictly hyperbolic systems of two conservation laws,
by Chen, et al. \cite{Chen,DCL} to fluid equations with  
$\gamma \in(1,5/3]$ and by Lions, Perthame, Souganidis, and Tadmor \cite{LPT,LPS} 
to the full range $\gamma >1$. 
The theory was extended to real fluid equations
by Chen and LeFloch \cite{CL1,LeFloch1,CL2}. We also mention the important 
work by Perthame and Tzavaras on the kinetic formulation for systems of two 
conservation laws; see \cite{Perthame,PT}.
The success of these works 
relies on a detailed analysis of the fundamental solution of the
entropy wave equation \eqref{2.2b}, which is a degenerate, linear wave equation.

When $\gamma =1$ the analysis developed in \cite{CL1,CL2} for the construction of entropies
does not work because the equation \eqref{2.2b} degenerates at a higher degree and
the Cauchy problem at the line $\rho =0$ becomes highly singular. One novelty of the present paper
is to rely on symmetry group argument to identify the entropy kernel.

For the convenience of the reader we summarize now the main steps of the proof of Theorems 2.1--2.3.

\

{\sl Step 1.} We rely on the vanishing viscosity method and first
construct a sequence of approximate solutions $(u^\eps,\rho^\eps)$, $\eps\downarrow 0$, defined on the strip
$\Pi$, and such that
$$
2 \, \eps^r \leq \rho^\eps \leq \rho_2 <1
$$
for some $r>1$.
The constant $\rho_2$ can be chosen to be arbitrary small by introducing a rescaled, initial density
$\lambda \, \rho_0$.
We will thus establish first Theorems 2.1 to 2.3 in the case when the initial density is small. Then we
will treat the general case by observing that the system (1.1)--(1.3)
is invariant via the symmetry
$(u,\rho) \to (u,\lambda\rho)$. More precisely, given an entropy solution $(u,\rho)$ of
the problem (1.1)--(1.3) with initial data $(u_0, \rho_0)$,
the pair $(u',\rho') := (u,\lambda\rho)$ is also an entropy solution with the initial data $(u'_0,\rho '_0)
= (u_0, \lambda\rho_0)$.

\

{\sl Step 2.} Next, we prove that there is a sequence $\eps\downarrow 0$ such that
\be
W^\eps := \rho^\eps e^{u^\eps} \rightharpoonup W,\quad
Z^\eps := \rho^\eps e^{-u^\eps} \rightharpoonup Z
\quad
\mbox{weakly $\star$ in } L_{loc}^\infty(\Pi)
\label{2.8}
\ee
and there exist Young measures $\nu_{x,t}$, associated with the sequence
$\eps\downarrow 0$ and defined on the $(W,Z)-$ plane for each point $(x,t)\in$
$\Pi$, such that
$$
\lim_{\eps\to 0}
F(W^{\eps }(x,t),Z^{\eps }(x,t)) = \iint F(\alpha, \beta)d\nu_{x,t}
=
\langle \nu_{x,t},F\rangle
=: \la F\ra\quad \mbox{ in } \DD'(\RR^2)
$$
for any $F(\alpha, \beta)\in C_{loc}(\RR^2)$.
The crucial point in the compensated compactness argument is to prove that $\nu$ is a point mass measure.
In that case the convergence in \eqref{2.8} becomes strong in any $L_{loc}^r(\Pi)$, $1\leq r<\infty$.

\

{\sl Step 3.} Given two entropy pairs $(\eta_i, q_i)$, obeying the conditions of Theorem 2.1, we check that
Tartar's commutation relations
\be
\la \nu_{x,t}, \eta_1\, q_2 - \eta_2 \,q_1 \ra \,
= \la  \nu_{x,t},\eta_1\ra\, \la \nu_{x,t},q_2\ra - \la \nu_{x,t},\eta_2\ra\,
\la \nu_{x,t},q_1\ra
\label{2.9}
\ee
hold. Here, we apply the so-called div-curl lemma of Murat \cite{Murat1} and Tartar
\cite{Tartar1,Tartar2}. The objective is to prove that
the measure $\nu$ is a point mass measure by using a ``sufficiently large'' class of entropy pairs in \eqref{2.9}.

\

{\sl Step 4.} To produce a large family of entropy pairs, we have to construct a fundamental solution
$\chi (R,u-s)$ (where $R:=\ln\rho$)
of the entropy equation \eqref{2.2b}. To this end, we rely on symmetry group arguments for
the equation \eqref{2.2b}. We find that it has an invariant solution
$$
\eta (u,\rho) = \sqrt{\rho} \, f(u^2 -\ln^2 \rho),\quad
\mbox{ where } \quad \xi \, f''(\xi) + f'(\xi) + \frac{1}{16} \, f(\xi) = 0.
$$
Then we define
$$
\chi = e^{R/2} \, f(|u-s|^2 - R^2) \, \mathbf{1}_{|u-s|< |R|}.
$$
The function $f(\xi)$ can be represented by a Bessel function of zero index.

\

{\sl Step 5.} Then we search for the entropy pairs in the form
\be
\eta = \int\chi (R,u-s) \, \psi (s) \, ds,
\qquad
q = \int\sigma (R,u,s) \, ds,
\label{2.10}
\ee
where $\psi\in L^1 (\RR)$ is arbitrary
and we describe properties of the kernels $\chi,\sigma$. In particular, we find that
$\sigma = u \, \chi (R,u-s)+ h(R,u-s)$, where the function $h$ is given by an explicit formula.
We also will show that
$$
P\chi :=\del_s\chi = e^{R/2} \, \Big( \delta_{s=u-|R|} - \delta_{s=u+|R|} \Big)
   + G^\chi(R,u-s) \, \1_{|u-s|<|R|},
$$
$$
Ph:= \del_s h = e^{R/2} \, \Big( \delta_{s=u-|R|} + \delta_{s=u+|R|} \Big)
       + G^h(R,u-s) \, \1_{|u-s|<|R|},
$$
where $G^\chi(R,v)$ and $G^h(R,v)$  are bounded, continuous functions.

\

{\sl Step 6.} Finally, we plug the entropy pairs \eqref{2.10} in Tartar's commutation relations,
but in the form derived by Chen and LeFloch \cite{CL1}. We arrive after cancellation of $\psi$
at the following equality in $\DD'(\RR)^3$
$$
\la\chi_1 P_2 h_2  -h_1 P_2\chi_2  \ra
\la P_3\chi_3 \ra +
\la h_1 P_3\chi_3 -  \chi_1 P_3 h_3 \ra\la P_2\chi_2 \ra
 = - \la P_3 h_3 P_2 \chi_2 - P_3\chi_3 P_2 h_2 \ra
\la\chi_1\ra,
$$
where the notations $g_i := g(R,u,s_i)$ and $P_i g_i := \del_{s_i}g(R,u,s_i)$ are used.
Then we test this equality with the function
$$
\frac{1}{\delta^2}\psi (s_1)\varphi_2 (\frac{s_1
-s_2}{\delta}) \varphi_3 (\frac{s_1 -s_3}{\delta}),
$$
where $\psi\in\DD(\RR)$ and $\varphi_j$ are molifiers such that
$$
Y:= \int_{-\infty}^\infty\int_{-\infty}^{s_2} \bigl(\varphi_2(s_2) \,
\varphi_3(s_3)-\varphi_3(s_2)\, \varphi_2(s_3)\bigr) \, ds_2ds_3 \ne 0.
$$
This identity involves products of measures by functions of bounded variation.
Such products were earlier discussed by Dal~Maso, LeFloch,
and Murat \cite{DLM}.

By letting $\delta$ go to zero we obtain the equalities
\be
Y\iint_{W,Z}D(\rho) \rho \iint_{\bigl\{W' < W\bigr\}
\cap \bigl\{Z'<1/W\bigr\}} \sqrt{\rho^{'}} \, d\nu (W',Z')d\nu
(W,Z)=0,
\label{2.11}
\ee
\be
Y\iint_{W,Z}D(\rho) \rho \iint_{ \bigl\{W' < 1/Z \bigr\}
\cap \bigl\{ Z'<Z \bigr\}} \sqrt{\rho^{'}} \, d\nu (W',Z') \,
d\nu (W,Z) = 0,
\label{2.12}
\ee
where
$$
\rho =(WZ)^{1/2},\quad
D(\rho)=\sqrt{\rho}(-\frac{1}{2}+\frac{15}{8}\ln\frac{1}{\rho}),\quad
\rho ' =(W'Z')^{1/2},
$$
and the measure $d\nu (W',Z')$ is a copy of $d\nu$ on the $(W',Z')$-plane.
At this point we choose the constant $\rho_2$ (see {\sl Step }1) small enough to ensure the
inequality $D(\rho)\geq \sqrt{\rho}/2$.
Hence, it follows from \eqref{2.11},\eqref{2.12} that
$d\nu_{x,t}=\alpha\delta_P +\mu_{x,t}$ and $\alpha \, (1-\alpha)=0$,
where $P(x,t)$ is a point on the $(W,Z)-$plane and the support of the measure $\supp \mu_{x,t}$
lies in the set$\{\rho =0\}$.
This representation formula for the measure $\nu_{x,t}$ enables us to justify the passage to
the limit as $\eps\downarrow 0$.
We summarize {\sl Step 6} in the following key result.

\begin{theorem}
Let $(m_n, \rho_n)$ be a bounded in $L^\infty(\Pi)$ sequence of entropy solutions of the problem
\eqref{2.1a} and such that
$$
0\leq\rho_n, \quad |m_n |\leq c\rho_n (1+|\ln\rho_n|)
$$
uniformly in $n$. Then, passing to a subsequence if necessary, $(m_n, \rho_n)$
converges almost everywhere in $\Pi$ to an entropy solution $(m, \rho)$ of \eqref{2.1a}.
\end{theorem}

\section{Vanishing viscosity method}
\setcounter{equation}{0}

Given parameters $\eps,\eps_1 >0$ we consider the Cauchy problem
\be
\rho_t +(\rho u)_x =\eps \rho_{xx}+2\eps_1 u_x,
\label{3.1}
\ee
\be
(\rho u)_t +(\rho u^2)_x + \rho_x =\eps (\rho u)_{xx}
+\eps_1 (u^2)_x +2\eps_1 (\ln \rho)_x,
\label{3.2}
\ee
with initial condition
\be
\rho |_{t=0}=\rho_0^\eps + 2 \, \eps_1, \qquad
u|_{t=0}=u_0^\eps.
\label{3.3}
\ee
In this section we establish the existence of smooth solutions to this problem.
Later in this section we will assume that $\eps_1 =\eps^r$ for some $r>1$.
The positivity of the density will be obtained by the following argument.

\begin{lemma} {\rm (Positivity for convection-diffusion equations.)}
If $v=v(x,t)$ is a
smooth bounded solution of the Cauchy problem
\be
v_t +(u \, v)_x =\eps v_{xx},
\qquad v|_{t=0}=v_0 (x),
\label{3.4}
\ee
where $u=u(x,t) \in L^\infty(\Pi)$ and $u_0 \in L^\infty(\RR)$,
then $v \geq 0$ provided $v_0 \geq 0$.
\end{lemma}

\begin{proof}  Given $R>0$, let $\psi:\RR_{+}\to\RR$ be a non-increasing function of class $C^2$
such that $\psi (x)=1$ for $x\in [0,R]$, $\psi (x)=e^{-x}$ for
$x\geq 2R$, and $\psi (x)$ is a non-negative polynomial for $R\leq
x\leq 2R$. Denote $\Psi (x)=\psi (|x|)$ for $x\in\RR$.
Clearly,
\be
|\Psi^{'}(x)|\leq\frac{c_1}{R}\Psi (x),\quad |\Psi
^{''}(x)|\leq\frac{c_1}{R^2}\Psi (x)
\label{3.5}
\ee
for some constant $c_1 >0$. The map
$$
U_{\mu} (v)=\left\{ \begin{array}{ll}
\sqrt {v^2 +\mu^2}-\mu, &  v\leq 0,\\
0, &  v>0,
\end{array} \right.
$$
is a regularization of the mapping $v \mapsto v_- := \max \{-v,0\}$.

Using \eqref{3.4} and \eqref{3.5} we can compute
 the $t$-derivative of the integral $\int\Psi U_{\mu}(v) \, dx$:
\be
\begin{split}
& \frac{d}{dt}\int\Psi U_{\mu} \, dx + \eps \, \int\Psi v_x^2 \frac{\del^2 U_{\mu}}{\del v^2} \, dx
\\
& = \int\frac{\del^2 U_{\mu}}{\del v^2} vv_x \, (\eps \, \Psi_x + u \Psi) \, dx
 + \int v\frac{\del U_{\mu}}{\del v }(\eps \Psi_{xx} + u \, \Psi_x) \, dx
\\
& \leq \int\frac{\del^2 U_{\mu}}{\del v^2} v |v_x
|\Psi (\eps c_1 /R +|u|)dx+ \int v |\frac{\del
U_{\mu}}{\del v }|\Psi (\eps c_1 /R^2+|u|c_1 /R) \, dx.
\end{split}
\label{3.6}
\ee
Observe that
$$
\eps v_x^2 -v |v_x|(\eps c_1 /R +|u|)=
\eps (|v_x|-v (\frac{c_1}{2R}
+\frac{|u|}{2\eps}))^2 - v^2
(\frac{c_1}{2R}+\frac{|u|}{2\eps})^2,
$$
$$
v^2\frac{\del^2 U_{\mu}}{\del v^2}\leq \frac{\mu
^2 v^2}{(v^2 +\mu^2)^{3/2}},
$$
and
$$
v\frac{\del U_{\mu}}{\del v }\to v_- \quad\mbox{as}\quad \mu\to 0.
$$
We integrate \eqref{3.6} with respect to $t$ and let $\mu$ tend to zero,
by taking into account that $U_{\mu}(v_0)=0$:
$$
\int\Psi v_- \, dx
\leq \int\limits_{0}^{t}\int\Psi v_-
(\frac{\eps \, c_1}{R^2}+ \frac{|u| \, c_1}{R}) \, dx d\tau.
$$
By Gronwall's lemma, $\int\Psi v_- \, dx=0$. We thus conclude that $v\geq 0$.
\end{proof}

As a consequence of Lemma 3.1, we deduce that any bounded solution $(u,\rho)$
of the problem \eqref{3.1}--\eqref{3.3} has the following property:
\be
\rho \geq 2 \, \eps_1 \quad \mbox{ uniformly in } \eps.
\label{3.7}
\ee
Namely, this is clear since the function $v=\rho -2 \, \eps_1$ solves the problem
$$
v_t +(u \, v)_x = \eps \, v_{xx}, \qquad v|_{t=0} \geq 0.
$$

From now on, we assume that the initial data $\rho_0^\eps$ and $u_0^\eps$ belong to the Sobolev space
$H^{2+\beta}(\RR)$ for some $0<\beta <1$ and satisfy
$$
0\leq \rho_0^\eps\leq M, \qquad \|u_0^\eps\|_\infty\leq u_1,
$$
and
$$
u_0^\eps\to u_0, \quad \rho_0^\eps \to \rho_0 \quad \mbox{ in } L_{loc}^{1}(\RR),
$$
where $u_1 : =\|u_0\|_\infty$ and $M :=\|\rho_0\|_\infty$.

\begin{lemma}
Let $(u,\rho)$ be a
smooth bounded solution of the Cauchy problem \eqref{3.1}--\eqref{3.3}. Then there
exist positive constants $c_1$, $\rho_1$, $W_1$, and $Z_1$ such that
$$
2 \, \eps_1 \leq \rho\leq\rho_1, \quad |m| := \rho \, |u| \leq c_1 \, \rho
\, (1+|\ln \rho|) \leq m_1,
\qquad \rho_1 := (2 \, \eps_1 + M) \, e^{u_1},
$$
\be
\begin{split}
& m_1 :=c_1 \sup_{0\leq\rho\leq\rho_1 } \rho
(1+|\ln \rho|),
\\
& 0\leq W:=\rho e^u \leq W_1, \quad 0\leq Z:=\rho e^{-u} \leq Z_1,
\end{split}
\label{3.8}
\ee
uniformly in $\eps$.
\end{lemma}

\begin{proof} Passing to the Riemann invariant variables
$$
w := u + \ln\rho, \quad z := u - \ln\rho,
$$
we can rewrite the system \eqref{3.1}-\eqref{3.2} as
$$
w_t + w_x \, \Big( u+1 -\frac{2\eps_1 }{\rho} +
\frac{\eps z_x}{2}-\frac{3\eps w_x}{4} \Big)
= \eps w_{xx}-\frac{\eps z_x^2}{4},
$$
$$
z_t +z_x (u-1 +\frac{2\eps_1 }{\rho} -
\frac{\eps w_x}{2}+\frac{3\eps z_x}{4
})= \eps z_{xx}+\frac{\eps w_x^2}{4} .
$$
By the maximum principle,
$$
w\leq \max w_0 (x),\quad z\geq \min z_0 (x).
$$
Now, the estimates (3.8) is a simple consequence of these
inequalities.
\end{proof}

By the estimates (3.8) there exist sequences $W^{\eps_n}$, $Z^{\eps_n}$,
 $\rho^{\eps_n}$, and
 $m^{\eps_n}:=$
$\rho^{\eps_n}u^{\eps_n}$
and  a family of non-negative probability measures
$\nu_{x,t}$, called Young measures, defined on the $(W,Z)$-plane, such that
\be
W^{\eps_n} \rightharpoonup W,\quad
Z^{\eps_n} \rightharpoonup Z,\quad
\rho^{\eps_n} \rightharpoonup \rho,\quad
\rho^{\eps_n}u^{\eps_n} \rightharpoonup m\quad
\mbox{ weakly $\star$ in } L_{loc}^\infty(\Pi),
\label{3.9}
\ee
and
$$
\iint_{\Pi} \Big( F(W^{\eps_n}(x,t),Z^{\eps_n}(x,t))
- \langle F\rangle\Big) \, \varphi (x,t) \, dtdx\to 0,
$$
where we have set $\langle F\rangle :=\int_{W,Z}F(W,Z)d\nu_{x,t}$
for any test function $\varphi \in\DD(\RR^2)$ and any continuous function
$F(W,Z)\in C_{loc}(\RR^2)$. Moreover,
$$
\mbox{supp}\, \nu_{x,t}\subset \{(W,Z):0\leq W\leq W_1, \quad
0\leq Z\leq Z_1
\}.
$$
For a proof that to each bounded sequence $v_n(x,t)$ one can associate 
a Young measure $\mu_{x,t}$, we refer to Tartar \cite{Tartar1} and 
Ball \cite{Ball}; see also \cite{Shelukhin1}.

\begin{lemma} {\rm (Entropy dissipation estimate.)}
The smooth solution
$(u,\rho)$ of the Cauchy problem \eqref{3.1}-\eqref{3.3} satisfies the estimate
\be
\|\frac{\eps\rho_x^2}{\rho}+\eps \, \rho u_x
^2\|_{L_{loc}^{1}(\Pi)}\leq c
\label{3.10}
\ee
uniformly in $\eps$.
\end{lemma}

\begin{proof} The identity
\be
\begin{split}
& \frac{\del}{\del t} \Big( \frac{\rho u^2 }{2}+ (1+\rho\ln \rho
-\rho) \Big)
+ \frac{\eps \rho_x^2}{\rho}+ \eps\rho u_x^2
\\
& =
- \frac{\del}{\del x} \{ \frac{\rho u^3}{2}+ u\rho\ln \rho
- \eps\rho_x \ln\rho - 2\eps_1  u\ln \rho -
\eps (\frac{\rho u^2}{2})_x -\frac{\eps_1 u^3}{3} \} =: - J_x
\end{split}
\ee
follows immediately from \eqref{3.1} and \eqref{3.2}. Multiplying this identity
by the function $\Psi (x)$ introduced in the proof of Lemma 3.1 and
integrating with respect to $x$ we deduce, in view of the estimates (3.7) and (3.8),
$$
\int J\Psi_x \, dx \leq \frac{1}{2}\int\Psi (\frac{\eps \rho_x
^2}{\rho}+ \eps\rho u_x^2)dx+c\int\Psi (1+\frac{\rho
u^2}{2})dx.
$$
Hence, we have
$$
\int_0^T \int\Psi (\frac{\eps \rho_x^2}{\rho}
+ \eps \, \rho u_x^2) \, dxdt \leq c,
$$
which yields the desired estimate.
\end{proof}

We rewrite the equations (3.1)-(3.2) as a quasi-linear parabolic system:
\be
u_t + a_1 (u,\rho, u_x \rho_x) = \eps \, u_{xx},
\quad \rho_t + a_2 (u,\rho, u_x \rho_x)=\eps \, \rho_{xx},
\label{3.11}
\ee
where we have set
$$
a_1 := uu_x -\frac{\rho_x}{\rho} - \frac{2\eps\rho_x
u_x}{\rho}- \frac{2\eps_1 \rho_x}{\rho^2},
\qquad
a_2 : =(\rho u)_x -2\eps_1 \, u_x .
$$
In view of (3.7) and (3.8), we obtain the global a priori estimates
$$
2 \, \eps_1 \leq \rho \leq \rho_1,
\qquad |u| \leq c(u_1,\rho_1,\eps_1).
$$

With these estimates at hand, it is a standard matter to derive estimates in H\"{o}lder's norms,
depending on $\eps$, by standard techniques of
the theory quasi-linear parabolic equations \cite{LSU}. We will only sketch the
derivation.  Let $\zeta (x,t)$ be a smooth function such that $\zeta \ne
0$ only if $x\in \omega$, where $\omega$ is an interval $[x_0 -\sigma, x_0 +\sigma]$.
Denote
$$
u^{(n)} := \max \{u-n,0\}.
$$
Multiplying the second equation in \eqref{3.11} by $\zeta^2 \rho^{(n)}$ and
integrating with respect to $x$, one obtains
$$
\frac{d}{dt}\int\zeta^2 |\rho^{(n)}|^2 dx+ \eps\int\zeta^2
|\rho_x^{(n)}|^2 dx\leq \gamma\int (\zeta_x^2 +\zeta |\zeta
_t|)|\rho^{(n)}|^2 dx+ \gamma \int\zeta \mathbf{1}_{\rho\geq n}dx.
$$
Similarly, for the velocity variable one gets
$$
\frac{d}{dt}\int\zeta^2 |u^{(n)}|^2 \, dx + \eps \, \int\zeta^2
|u_x^{(n)}|^2 \, dx \leq \gamma
\int (\zeta_x^2 +\zeta |\zeta
_t|)|u^{(n)}|^2 + \zeta \mathbf{1}_{\rho\geq n}+
\eps\zeta^2 |\rho_x^{(n)}|^2 \, dx.
$$
These inequalities imply  that $u$ and $\rho$ belong to a class
$\mathcal{B}_2 (Q,M,\gamma, r,\delta, n)$ \cite{LSU} (Chapter II,
\S 7, formula $(7.5)$), for some parameters $Q, M, \gamma,  r, \delta $,
and $n$.
Then it follows that the estimate
$$
\|u,\rho\|_{H^{\alpha, \alpha /2}(\omega\times [0,T])}\leq c
$$
holds for some $\alpha\in (0,1)$.

In the same manner, one can estimate the H\"{o}lder norm of the
derivatives $u_x$, $u_{xx}$, $u_t$, $\rho_x$, $\rho_{xx}$, and
$\rho_t$, in the same way as done in \cite{FS}
for a general class of parabolic systems.

We now arrive at the main existence result, concerning the viscous approximation (3.1)-(3.3).

\begin{lemma} {\rm (Existence of smooth solution of the regularized system.)}
Let $u_0
^\eps$, $\rho_0^\eps$ $\in L^\infty\cap
H_{loc}^{\beta}$, $0<\beta <1$. Then the Cauchy problem (3.1)-(3.3)
has a unique solution such that
$$
u,\rho \in L^\infty(\Pi)\cap H_{loc}^{2+\beta
,1+\beta /2} (\Pi).
$$
\end{lemma}
Now, we set $\eps_1 =\eps^r$, $r>1$, and study compactness
of the viscous solutions $(u^\eps,\rho^\eps)$
when $\eps \to 0$.

\begin{lemma}
Given an entropy
entropy-flux pair $(\eta (m,\rho),q(m,\rho))$, $m=\rho u$, the
sequence
$$
\theta^\eps :=  \frac{\del \eta^\eps}{\del t}+\frac{\del
q^\eps}{\del x}
$$
is compact in $W_{loc}^{-1,2}(\Pi)$, where $
\eta^\eps=\eta (m^\eps,\rho^\eps)$, $
q^\eps=q(m^\eps,\rho^\eps). $
\end{lemma}

\begin{proof}
We use the following lemma due to Murat's lemma \cite{Murat2}.

\vspace{0.3cm}
\noindent
{\it
Let $Q\subset\RR^2$ be a bounded domain, $Q\in C^{1,1}$.
Let $A$ be a compact set in
$W^{-1,2}(Q)$, $B$ be a bounded set in the space of bounded Radon measures $M(Q)$, and $C$ be a bounded
set in $W^{-1,p}(Q)$ for some $p\in (2,\infty]$.
Further, let $D\subset\DD'(Q)$ be such that
$$
D\subset (A+B)\cap C.
$$
Then there exists  $E$, a compact set in $W^{-1,2}(Q)$ such that
$D\subset E$.
}

By definition, the functions
$\eta (m,\rho)$ and $q(m,\rho)$ solve the system
$$
q_m =\frac{2m}{\rho}\eta_m  +\eta_{\rho},\quad q_{\rho}=\eta_m
-\frac{m^2}{\rho^2}\eta_m  .
$$
Hence, calculations show that
\be
\theta^\eps= 2\eps_1 m_x (\frac{\eta
_{\rho}^\eps}{\rho}+\frac{m\eta_m^\eps}{\rho
^2})+ 2\eps_1 (-\frac{m\eta_{\rho}^\eps}{\rho
^2}- \frac{m^2\eta_{m}^\eps}{\rho^3}+\frac{\eta
_{m}^\eps}{\rho}) +\eps\eta
_{\rho}^\eps\rho_{xx}+\eps\eta_m^\eps \, m_{xx}
=
\label{3.12}
\ee
$$
\eps_1 u_x (q_m^\eps+\eta
_{\rho}^\eps)- 2\eps_1 \frac{\rho_x \eta_m
^\eps}{\rho}+ \eps\eta_{xx}^\eps-
\eps [\eta_{\rho \rho}^\eps\rho_x^2 +\eta_{mm}^\eps \, m_x^2 +2\eta_{\rho m}^\eps\rho_x
m_x ].
$$
We check the conditions of Murat's lemma.
By  Lemma 3.2, the sequence $\theta^\eps$ is bounded in
$W_{loc}^{-1,\infty}(\Pi)$. Hence, it is enough to show that $\eps \eta_x
^\eps$ $\to 0$ in $L_{loc}^2(\Pi)$
and the residual sequence $\theta^\eps-\eps \eta
_{xx}^\eps$ is bounded in $L_{loc}^{1}(\Pi)$.

We have
$$
\eps\eta_x^\eps= \eps\rho u_x \eta_m
^\eps+ \eps \rho_x \frac{q_m^\eps+\eta
_{\rho}^\eps}{2}.
$$
Thus, by estimates (3.8) and (3.10), $\eps \eta_x
^\eps\to 0$ in $L_{loc}^2$.

Consider the sequence $\theta^\eps-\eps \eta
_{xx}^\eps$. We have
$$
\theta^\eps-\eps \eta_{xx}^\eps=-
\eps [\eta_{\rho \rho}^\eps\rho_x^2 +\eta
_{mm}^\eps m_x^2 +2\eta_{\rho m}^\eps\rho_x
m_x ]+ \eps_1 u_x (q_m^\eps+\eta
_{\rho}^\eps) -2\eps_1 \frac{\rho_x \eta_m
^\eps}{\rho} .
$$
Each term on the right hand-side is bounded in $L_{loc}^{1}$ provided
$\eps_1 =\eps$. Indeed, by (3.7),
$$
2\eps_1 |u_x| =\frac{2\eps_1 \rho
^{1/2}|u_x|}{\rho^{1/2}}\leq \sqrt{2\eps}\rho^{1/2}|u_x
|,\quad \frac{2\eps_1 |\rho
_x|}{\rho}\leq\frac{\sqrt{2\eps}|\rho_x|}{\rho^{1/2}} .
$$
Moreover, if $\eps_1 =0(\eps)$,
\be
\eps_1 u_x (q_m^\eps+\eta_{\rho}^\eps)
-2\eps_1 \frac{\rho_x \eta_m^\eps}{\rho}\to
0\quad \mbox{in}\quad L_{loc}^2(\Pi).
\label{3.13}
\ee
The other terms are treated similarly. This completes the proof.
\end{proof}

Given two entropy pairs $(\eta_i (m,\rho),q_i (m,\rho))$, $(i=1,2)$, from Lemma 3.5,
we define
$$
\tilde{\eta}_i (W,Z)=\eta_i (f_2 (W,Z),f_1 (W,Z)),\quad \tilde{q}_i (W,Z)=q_i (f_2 (W,Z),f_1 (W,Z)).
$$
Clearly, the functions
$$
\del_t \tilde{\eta}_i^\eps + \del_x \tilde{q}_i^\eps
$$
are compact in $W_{loc}^{-1,2}(\Pi)$.
Hence, by the div-curl lemma \cite{Tartar1}, Tartar's commutation relation
\be
\la \tilde{ \eta}_1\, \tilde{q}_2 - \tilde{\eta}_2 \,\tilde{q}_1 \ra \,
= \la  \tilde{\eta}_1\ra\, \la \tilde{q}_2\ra - \la \tilde{\eta}_2\ra\, \la \tilde{q}_1\ra
\label{3.14}
\ee
is valid.

For reader's convenience, we remind that the div-curl lemma states the following.

\vspace{0.2cm} \noindent
{\it Let $Q\subset\RR^2$
be a bounded domain, $Q\in C^{1,1}$. Let
$$
w_1^k \rightharpoonup w,\quad
w_2^k \rightharpoonup w_2,\quad
v_1^k \rightharpoonup v_1,\quad
v_2^k \rightharpoonup v_2,
$$
weakly in $L^2 (Q)$, as $k\to\infty$.
With
$\mbox{curl}(w_1, w_2)$ denoting
$\del w_2 /\del x_1 -$
$\del w_1 /\del x_2$,
suppose that
the sequences
$\mbox{div}(v_1^k, v_2^k)$
and $\mbox{curl}(w_1^k, w_2^k)$
lie in a compact subset $E$ of $W^{-1,2}(Q)$.
Then, for a subsequence,
$$
v_1^k w_1^k +v_2^k w_2^k \to v_1  w_1 +v_2  w_2 \quad
\mbox{in $\DD'(Q)$}\quad
\mbox{as}\quad k\to\infty.
$$
}

The further claim is due to the fact that system
\eqref{1.1}-\eqref{1.3}  is
invariant with respect to the scaling $\rho\to\lambda\rho$.

\vspace{0.2cm} \noindent
\begin{lemma}
If $(m,\rho)$ is an
entropy solution with initial data $(m_0, \rho_0)$ then
$(cm,c\rho)$ is also the entropy solution with the initial data
$(cm_0, c\rho_0)$, where $c$ is an arbitrary positive constant.
\end{lemma}
\vspace{0.2cm} \noindent {\it Proof.}
The claim follows easily
from the fact that the pair $(\eta (cm,c\rho),q(cm,c\rho))$ is an
entropy-entropy flux pair as soon as the pair $(\eta
(m,\rho),q(m,\rho))$ is an entropy-entropy flux pair.

Given $\lambda >0$, let us consider the auxiliary problem
\be
\rho_t +(\rho u)_x =\eps \rho_{xx}+2\eps_2 u_x,
\label{3.15}
\ee
\be
(\rho u)_t +(\rho u^2)_x +\rho_x =\eps (\rho u)_{xx}
+\eps_2 (u^2)_x +2\eps_2 (\ln \rho)_x,
\label{3.16}
\ee
\be
\rho |_{t=0}=\lambda\rho_0^\eps(x)+2\eps_2, \quad
u|_{t=0}=u_0^\eps(x),
\label{3.17}
\ee
where $\eps_2 =\lambda\eps_1=\lambda\eps^r$.

The main feature of the auxiliary problem is the following. If the functions $(u_\eps,\rho_\eps)$
solve the problem (3.1)-(3.3) then the functions
$(u_\eps,\rho_\eps^{'})$ solve the problem \eqref{3.15}-\eqref{3.17} with $\rho_\eps^{'}=$
$\lambda\rho_\eps$.

The solution
$(u_\eps,\rho_\eps)$
of problem \eqref{3.15}-\eqref{3.17}  obeys the estimates
\be
2\eps_2 \leq \rho_\eps\leq
(2\eps_2 +\lambda\|\rho_0 \|_\infty)e^{\|u_0\|_\infty} =: \rho_2,
\quad
|u_\eps\rho_\eps|\leq c\rho_\eps(1+|\ln \rho_\eps|)
\label{3.18}
\ee
uniformly in $\eps$.
Lemmas 3.3-3.5 are also valid for $(u_\eps,\rho_\eps)$.
The corresponding Young measure $\nu_{x,t}$ has a finite support:
\be
\supp \nu_{x,t}\subset \{(W,Z):0\leq W\leq W_2, \quad
0\leq Z\leq Z_2 \} := K.
\label{3.19}
\ee
We impose the following smallness conditions for $\lambda$:
\be
\rho_2 <1,\quad
\ln \frac{1}{\rho_2}\geq \frac{8}{15}.
\label{3.20}
\ee
Assume that the solution $(u_\eps,\rho_\eps)$ of the auxiliary problem converges
to an entropy solution $(m,\rho)$ of the problem \eqref{2.1a}:
$$
(u_\eps\rho_\eps,\rho_\eps)\to (m,\rho)\quad\mbox{almost everywhere in}\quad \Pi .
$$
The initial data for $(m,\rho)$ are
$$
\rho |_{t=0}=\lambda\rho_0,\quad
m|_{t=0}=\lambda m_0 .
$$
By  Lemma 3.6, the functions $(m^{'},\rho^{'})=$
$(m/\lambda, \rho /\lambda)$ is an entropy solution of the same problem with the initial data
$$
\rho^{'}|_{t=0}=\rho_0,\quad
m^{'}|_{t=0}= m_0 .
$$
Thus, it is enough to study convergence of the solutions to the auxiliary problem.

With the condition \eqref{3.20} at hand, the function
$$
D(R):=
(-\frac{1}{2}+\frac{15|R|}{8})e^{R/2},\quad R:=\ln \rho,
$$
from  Section 5 admits the estimate
$
D(R)\geq\frac{1}{2} e^{R/2}.
$
Hence, $D(R)$ vanishes only at the vacuum points $\rho =0$.

To conclude the section, we remark that the parameter $\eps_1$ serves as a regularizer
for the hyperbolic system \eqref{1.1}-\eqref{1.3} with
$\eps =0$ due to the estimate \eqref{3.7} (cf. \cite{Lu}).


\section{A large class of mathematical entropies}
\setcounter{equation}{0}

\subsection{Symmetry group analysis}

We already pointed out that a pair $(\eta,q)$ is a mathematical entropy if and only if
$\eta$ satisfies
\be
\eta_{\rho \rho} = \frac{1}{\rho^2}\eta_{uu}.
\label{4.1}
\ee
In order to derive an explicit formula for the weak entropies of the Euler system
we rely on symmetry group analysis, following \cite{Shelukhin2}.
Using the Riemann invariants
$$
w := u + \ln{\rho}, \quad z := u - \ln \rho,
$$
we write the equation \eqref{4.1} for the entropies in the form
\be
F(\eta_w, \eta_z, \eta_{wz})
:= \eta_{wz} - A \, (\eta_z - \eta_w) = 0,
\quad
A : = \frac{1}{4}.
\label{7.1}
\ee
In the more general case where $A$ is a function of $w$ and $z$, complete group analysis
arguments were developed in Ovsyannikov's monograph \cite{Ovsyannikov}.
In our case, $A$ is a constant and this analysis is simpler.
We only present the results of the formal derivation and we refer to \cite{Ovsyannikov}
for further details on the theory.

We look for a one-parameter group determined by the infinitesimal operator
$$
X=\xi (w,z,\eta)\frac{\del}{\del w}+\tau (w,z,\eta)\frac{\del}{\del z}+\varphi (w,z,\eta)
\frac{\del}{\del \eta}.
$$
Calculation of the first and the second prolongations of this operator yields
$$
X^1 =X+\zeta^{\eta_w}\frac{\del}{\del \eta_w}+
\zeta^{\eta_z}\frac{\del}{\del \eta_z},\quad
X^2 =X^1 +\zeta^{\eta_{ww}}\frac{\del}{\del \eta_{ww}}+
\zeta^{\eta_{wz}}\frac{\del}{\del \eta_{wz}}+
\zeta^{\eta_{zz}}\frac{\del}{\del \eta_{zz}},
$$
where
$$
\zeta^{\eta_{w}}=
D_w \varphi -\eta_w D_w \xi -\eta_z D_w \tau,\quad
D_w =\frac{\del}{\del w}+\eta_w \frac{\del}{\del \eta},
$$
$$
\zeta^{\eta_{z}}=
D_z \varphi -\eta_w D_z \xi -\eta_z D_z \tau,\quad
D_z =\frac{\del}{\del z}+\eta_z \frac{\del}{\del \eta},
$$
and
\be
\begin{split}
\zeta^{\eta_{wz}}
=
& D_z \varphi_w +\eta_w D_z \varphi_{\eta} +\varphi_{\eta} \eta_{wz}-
\eta_{ww}D_z \xi -\eta_{wz} (D_w \xi+D_z \tau)
\\
&
-\eta_{w}(D_z \xi_{w}+\eta_w D_z \xi_{\eta}+\xi_{\eta }\eta_{wz})
-\eta_{zz}D_w \tau
-\eta_{z}(D_z \tau_{w}+\eta_w D_z \tau_{\eta}+\eta_{wz }\tau_{\eta}).
\end{split}
\nonumber
\ee
Note that we need not calculate the coefficients
$\zeta^{\eta_{ww}}$ and $\zeta^{\eta_{zz}}$.
Application of the operator $X^2$ to $F$ and analysis of this application on the manifold $F=0$
enable us to conclude that the equation \eqref{7.1} admits four one-dimensional groups $G_i$
and one infinite-dimensional group $G_5$, associated with the infinitesimal operators
$$
\frac{\del}{\del w},\quad
\frac{\del}{\del z},\quad
\eta\frac{\del}{\del \eta},\quad
E :=  w\frac{\del}{\del w}-
z\frac{\del}{\del z}+
A(w+z)\eta\frac{\del}{\del \eta},\quad
\beta (w,z)\frac{\del}{\del \eta},
$$
where $\beta$ is a solution to \eqref{7.1}. The fact that
the equation \eqref{7.1} admits the group $G_i$
means the following:
if $\eta (w,z)$ solves \eqref{7.1} then for any $c,\xi\in R$
the following functions are solutions of \eqref{7.1} as well:
$$
\eta (w+c,z),\quad \eta (w,z+c),\quad c\eta (w,z),
$$
$$
\eta (e^{-\xi}w,e^{\xi}z)\exp {(Aw(1-e^{-\xi})-Az(1-e^{\xi}))},\quad
\eta (w,z)+\beta (w,z).
$$
Note that, once this assertion is obtained, its validity can be checked directly without 
referring to group analysis.

Let us find an invariant solution to the equation \eqref{7.1},
by applying the one-dimensional group $G_4$ associated with the infinitesimal operator
$E$.
First, we look for invariants $I(w,z,\eta)$ of the group $G_4$ as solutions of the equation
$E(I) = 0$. By the method of characteristics, one derives the system of O.D.E.'s
$$
\frac{dw}{w}=-\frac{dz}{z}=\frac{d\eta}{A(w+z)}
$$
and obtains easily the following two invariants:
$$
I_1 = wz, \quad I_2 =\eta e^{-A(w-z)}.
$$
Then, we look for an invariant solution of equation \eqref{7.1} in the form
(see again \cite{Ovsyannikov})
$I_2 = f(I_1)$,
or equivalently
\be
\eta =e^{A(w-z)}f(wz) .
\label{7.2}
\ee
Plugging \eqref{7.2} in \eqref{7.1}, we arrive at the following condition for the function $f(x)$:
\be
xf''(x) + f'(x) + A^2 f(x) = 0.
\label{7.3}
\ee
In conclusion, the equation \eqref{7.1} admits the solution
$$
\eta =\rho^{1/2} \, f(u^2 - \ln^2 \rho),
$$
where the function $f$ satisfies the equation \eqref{7.3}.


\subsection{Mathematical entropies}

We search for entropies $\eta = \eta(\rho,u)$ having the form
$$
\eta(\rho,u) = \rho^{1/2} \, f(u^2 -  \ln^2 \rho). 
$$
It is straightforward to see that $\eta$ solves the entropy equation \eqref{4.1}
if and only if the function $f=f(m)$ is a solution of the ODE
\be
m f''+f'+A^2 f=0,\quad A^2 = \frac{1}{16}.
\label{4.2}
\ee
With the notation
$$
R := \ln \rho
$$
the entropy then takes the form
$$
\eta = \eta(R,u) = e^{R/2} \, f(u^2 - R^2),
$$
while the entropy equation \eqref{4.1} reads
\be
\LLL (\eta) := \eta_{RR} - \eta_{uu} - \eta_{R} =0.
\label{4.3}
\ee

One solution to the second-order equation \eqref{4.2} is given by the following
expansion series
$$
f(m) := \sum_{n=0}^\infty \bigl( {A^n \over n!} \bigr)^2 (-1)^n m^n,
$$
with
$$
f(0)=1, \quad f'(0) = - A^2, \qquad
f(-y^2) = \sum_{n=0}^\infty(\frac{A^n y^n}{n!})^2.
$$
Observe that $f(m)$ can be represented by the Bessel function of zero order.
Given any function $g: \RR \to \RR$, we introduce the notation
$$
\overline g(m) = \begin{cases}
g(m),        &  m \leq 0,
\\
0,           &  m > 0.
\end{cases}
$$
In particular, this defines the function $\overf$. We denote by $\delta$ the Dirac
measure concentrated at the point $m=0$ and, more generally, by $\delta_{m=a}$ the Dirac measure
concentrated at the point $a$. We denote by $\DD(\RR)$ the space of smooth functions
with compact support and by $\DD'(\RR)$ the space of distributions.

\begin{lemma}
The function $\overf$ solves the ordinary differential equation \eqref{4.2} in $\DD'(\RR)$.
\end{lemma}

\begin{proof}
Given a test function $\varphi \in \DD(\RR)$, we have
$$
\la m \, \overf'', \varphi \ra
 := \int_\RR (m\varphi)'' \, \overf \, d m = \int_{-\infty}^0 (m\varphi)'' \, f \, d m
 = \la f(0) \, \delta + \overline{m \, f''}, \varphi\ra
$$
and
$$
\la\overf', \varphi\ra =
\la -f(0) \, \delta + \overline{f'}, \varphi\ra.
$$
Hence, we find
$$
\la m \, \overf'' + \overf' +A^2 \, \overf, \varphi \ra
=
\la \overline{m \, f'' + f' + A^2 \, f}, \varphi \ra
= 0.
$$
\end{proof}

Motivated by Lemma 4.1 we introduce the function
\be
\begin{split}
\chi(R,u)
:= e^{R/2} \, \overf(u^2 -R^2)
& = e^{R/2} \, \1_{R^2 -u^2 \geq 0} \, f(R^2 -u^2)
\\
& = e^{R/2} \, \sum_{n=0}^\infty \bigl( \frac{A^n }{n!} \bigr)^2 \, (R^2 - u^2)_+^n,
\end{split}
\label{4.4}
\ee
where
$$
\lambda_+ :=
\begin{cases}
\lambda,      & \lambda \geq 0,
\\
0,            & \lambda < 0,
\end{cases}
$$
and $\1_{g \geq 0}$ denotes the characteristic function of the set $\big\{ g \geq 0\big\}$.

\begin{theorem} {\rm (Existence of the entropy kernl.)}
The function $\chi$ defined by \eqref{4.4} is a
fundamental solution of the equation \eqref{4.3} in $\DD'(\RR^2)$. More precisely,
we have
$
\LLL (\chi) = 4 \, \delta_{(R,u)=(0,0)}.
$
\end{theorem}

\begin{proof} From the definition \eqref{4.4} of $\chi$
and since the multiplicative factor $e^{R/2}$
is a smooth function, it is straightforward to obtain
$$
\LLL(\chi) = e^{R/2} \, \bigl( \overf_{RR} - \overf_{uu} - {\overf \over 4} \bigr)
$$
in the sense of distributions. Note that, throughout the calculation,
$
f = f(u^2 - R^2)
$
and that the variables $u$ and $R$ describe $\RR$. 
We compute each term in the right-hand side of the above identity successively.
We have first
\be
\begin{split}
\la\overf_{RR}, \varphi\ra & = \la\overf, \varphi_{RR}\ra
= \iint_{ u^2 - R^2 \leq 0} f \, \varphi_{RR}\, dudR
\\
& = \iint_{|R|>|u|} f \, \varphi_{RR} \, dRdu
\\
& = \iint_{|R|>|u|} (\varphi_R \, f)_R + 2R \, f' \, \varphi_R \, dRdu .
\end{split}
\nonumber
\ee
Hence, we obtain
\be
\begin{split}
& \la\overf_{RR}, \varphi\ra
\\
& = f(0) \, \int_\RR \bigl( \varphi_R (-|u|,u) - \varphi_R (|u|,u) \bigr) \, du
+
\iint_{|R|>|u|} \bigl( (2R \, \varphi f')_R - \varphi \, (2 f' - 4R^2 \, f'') \bigr) \, dRdu
\\
& = f(0) \,  \int_\RR \bigl( \varphi_R (-|u|,u) - \varphi_R (|u|,u) \bigr) \, du
- 2 f'(0) \,  \int_\RR \bigl( \varphi(-|u|,u) \, |u| + \varphi(|u|,u) \, |u| \bigr) \, du
\\
& \hskip8.cm + \iint_{ |R|>|u| } \varphi \, (4 R^2 \, f'' - 2 f') \, dRdu.
\end{split}
\nonumber
\ee
Thus, we have established that
\be
\overf_{RR} = \overline{4R^2 \, f'' - 2 \, f'} + J_1,
\label{4.5}
\ee
where $J_1$ is the distribution defined by
\be
\begin{split}
\la J_1, \varphi\ra =
& f(0) \, \Big( \int_{-\infty}^0 \bigl( \varphi_R (u,u) - \varphi_R (-u,u) \bigr) \, du
+ \int_0^{+\infty} \bigl( \varphi_R (-u,u) - \varphi_R (u,u) \bigr) \, du \Big)
\\
&
+  2 \, f'(0) \, \Big(\int_{-\infty}^0 \bigl( u \, \varphi (u,u)
            + u \, \varphi(-u,u) \bigr) \, du
  - \int_0^{+\infty} \bigl( u \, \varphi(-u,u) + u \, \varphi(u,u) \bigr) \, du \Big).
\end{split}
\nonumber
\ee

The derivative in $u$ is determined in a completely similar fashion. We get
\be
\begin{split}
\la \overf_{uu}, \varphi \ra
= & \iint_{|u| < |R|} f \, \varphi_{uu} \, du dR
\\
= & \int_{-\infty}^{+\infty}\int_{-|R|}^{+|R|} \bigl( (\varphi_u f)_u - 2 u \, \varphi_u \, f' \bigr) \, dudR
\\
= & f(0) \, \int\limits_{-\infty}^{+\infty} \bigl( \varphi_u (R,|R|) - \varphi_u(R,-|R|) \bigr) \, dR
\\
& + \iint\limits_{|u| < |R|} \bigl( \varphi \, ( 2 f' + 4 u^2 \, f'')
      - 2 \, (u \, f' \, \varphi)_u \bigr) \, dudR
\end{split}
\nonumber
\ee
and thus
\be
\overf_{uu} = \overline{2 f' + 4 u^2 \, f''} + J_2,
\label{4.6}
\ee
where the distribution $J_2$ is given by 
$$
\la J_2, \varphi\ra
= f(0) \, \Big( \int_{-\infty}^0\varphi_u (R,-R)-\varphi_u (R,R)dR+
\int_0^{+\infty}\varphi_u (R,R)-\varphi_u (R,-R)dR
\Big)
$$
$$
+ 2 f'(0) \, \Big( \int_{-\infty}^0 \bigl( R \, \varphi (R,-R) + R\, \varphi (R,R) \bigr) \, dR
  - \int_0^{+\infty} \bigl( R \, \varphi(R,R) + R \, \varphi(R,-R) \bigr) \, dR \Big).
$$

Now, since the function $f$ satisfies the differential equation \eqref{4.2}
it follows from \eqref{4.5}-\eqref{4.6} that
$$
\overf_{RR} - \overf_{uu} - \overf/4 = J_1 - J_2.
$$
To conclude, we observe that
$$
\int_{-\infty}^0 u \, \varphi (-u,u) \, du = - \int_0^\infty R \, \varphi (R,-R) \, dR,
$$
$$
\int_0^\infty u \, \varphi (-u,u) \, du = - \int_{-\infty}^0 R \, \varphi (R,-R)dR,
$$
$$
\frac{d}{dR}\varphi (R,R) = \varphi_u (R,R) + \varphi_R(R,R),
\qquad
\frac{d}{dR}\varphi (R,-R) = -\varphi_u (R,-R) + \varphi_R (R,-R),
$$
and
$$
\int_{-\infty}^0 \varphi_u (R,R)\, dR = \varphi (0) - \int_{-\infty}^0 \varphi_R (R,R) \, dR,
$$
$$
\int_0^\infty\varphi_u (R,R) \, dR = - \varphi (0) - \int_0^\infty \varphi_R (R,R) \, dR,
$$
$$
\int_{-\infty}^0 \varphi_u (R,-R) \, dR = -\varphi (0) + \int_{-\infty}^0 \varphi_R (R,-R) \, dR,
$$
$$
\int_0^\infty\varphi_u (R,-R) \, dR = \varphi (0) + \int_0^\infty \varphi_R (R,-R) \, dR,
$$
we find that
$$
\la J_1 - J_2, \varphi \ra = 4 f(0) \, \varphi (0).
$$
Since $f(0) = 1$ and $e^{R/2}=1$ when $R=0$, this completes the proof of Lemma 4.2.
\end{proof}

\


\begin{lemma} The kernel $\chi$ vanishes on the vacuum
$$
\lim_{R \to -\infty} \chi (R,u) = 0 \mbox{ for every } u, 
$$
and, at the origin $R=0$, satisfies
\be
\lim_{R \to 0} \chi (R,\cdot) = 0,
\qquad
\lim_{R \to 0 \pm} \chi_R (R,\cdot) \to \pm 2 \, \delta_{u=0}
\label{4.7}
\ee
in the distributional sense in $u$.
Moreover, for any fixed $R$, $\chi$ has a compact support, precisely
$$
\chi (R,u) = 0, \qquad |u| > R.
$$
It is smooth everywhere except along the boundary of its support
where it has a jump of strength $\pm e^{R/2}$.
\end{lemma}

\begin{proof}
Detailed behavior of $\chi$  as $R\to -\infty$ can be derived from the asymptotic formula
\cite{Olver}
$$
\sum\limits_{0}^\infty(\frac{x^n}{n!})^2 =
\frac{e^{2x}}{2\sqrt{\pi x}}(1+O(\frac{1}{x})) \quad \mbox{ when } x \uparrow \infty.
$$
It follows that
$$
\chi (R,u)=
\1_{R^2 -u^2 \geq 0}
\frac{
e^{(-|R|+\sqrt{R^2 -u^2})/2}}
{\sqrt{\pi }(R^2 -u^2)^{1/4}} \,
(1+O(\frac{1}{\sqrt{R^2 -u^2}})) \quad \mbox{ when } R \downarrow -\infty.
$$
Next, given $\varphi=\varphi(u), \psi=\psi(R) \in\DD(\RR)$ we have
\be
\begin{split}
\la \chi_R, \varphi \, \psi \ra
& =
-\int_\RR \varphi (u) \, \int_{|R|>|u|} e^{R/2} \, f(m) \, \psi_R \, dRdu
\\
& = - J + \int_{|R|>|u|} \varphi \, \psi \, e^{R/2} \, \bigl( \frac{f}2 - 2 R \, f' \bigr) \, dudR,
\end{split}
\nonumber
\ee
where
\be
\begin{split}
J
& = \int_\RR \varphi (u) \Big\{ \int_{-\infty}^{-|u|} + \int_{|u|}^\infty \Big \} \, ( e^{R/2} \, f\psi)_R \, dRdu  \\
& = \int_\RR \varphi (u) \, \Big( e^{-|u|/2} \, \psi (-|u|) - e^{|u|/2} \, \psi (|u|) \Big) \, du.
\end{split}
\nonumber
\ee
We calculate
$$
\int_\RR\varphi (u) \, e^{-|u|/2} \, \psi (-|u|) \, du
= \int_\RR e^{R/2}\psi (R) \, \bigl( \varphi (R) + \varphi (-R) \bigr) \, \1_{R<0} \, dR,
$$
$$
\int_\RR \varphi (u) \, e^{|u|/2} \, \psi (|u|) \, du
= \int_\RR e^{R/2} \psi (R) \, \bigl( \varphi (R) + \varphi (-R) \bigr) \, \1_{R>0} \, dR.
$$
It follows that, for each $R$, $\chi_R$ is a distribution in the variable $u$,
given by the formula
\be
\begin{split}
\la \chi_R (R,\cdot), \varphi (u) \ra
= &
\int_{|u|<|R|} \varphi (u) \, e^{R/2} \, \Big(\frac{f(u^2 - R^2)}2 - 2 R \, f'(u^2 -R^2) \Big) \, du
\\ & + e^{R/2} \bigl( \varphi (R) + \varphi(-R) \bigr) \, \bigl( \1_{R>0} - \1_{R<0} \bigr).
\end{split}
\nonumber
\ee
This completes the proof of \eqref{4.7} and thus the proof of Lemma 4.3.
\end{proof}

Since the equation \eqref{4.3} is invariant under the transformations
$u \mapsto u - s$ for every constant $s$, we deduce immediately
from Lemma 4.2 that, for every $s \in \RR$, the function
$$
\chi(R,u-s) = e^{R/2} \, \overf(|u-s|^2 -R^2)
$$
satisfies the partial differential equation
\be
\LLL(\chi)(R,u-s) = 4 \, \delta_{(R,u)=(0,s)}
\label{4.8}
\ee
in $\DD'(\RR^2)$. We arrive at :

\begin{theorem} {\rm (The class of weak entropies to the isothermal Euler equations.)}
Restrict attention to the region $R < 0$ (respectively, $R>0$). The formula
$$
\eta(R,u) = \int_\RR \chi(R,u-s) \, \psi(s) \, ds,
$$
where $\psi$ is an arbitrary function in $L^1(\RR)$ describes the class of all
weak entropies to the Euler equations for isothermal fluids $(1.1)$--$(1.3)$.
In particular, for all $u \in \RR$ we have
\be
\lim_{R\to 0} \eta(R,u) = 0,
\quad
\lim_{R \to 0 \pm} \eta_R (R,u) = \pm 2 \, \psi(u),\quad
\lim_{R\to -\infty} \eta(R,u) = 0.
\label{4.9}
\ee
\end{theorem}

\begin{proof}
It follows from \eqref{4.8} that, for all $\varphi \in \DD(\RR^2)$,
$$
\int_\RR \LLL(\eta) \varphi \, dRdu = 4 \, \int_\RR \psi(s) \, \varphi (s,0) \, ds,
$$
which implies that
$$
\LLL(\eta) =0, \quad R \neq 0.
$$
Since, for any fixed $s, R$, the fundamental solution
$\chi (R,u-s)$ has a compact support in the variable $u$, we also have
$$
\int_\RR \chi(R,u-s) \, \psi (s) \, ds
\to 0, \qquad  R \to 0.
$$
\end{proof}


\subsection{Mathematical entropy-flux functions}

We look for the entropy-flux kernel $\sigma$ which should generate the class of
entropy flux-functions $q$ via the general formula
$$
q(R,u) = \int_\RR \sigma (R,u,s) \, \psi(s) \, ds.
$$
In the variables $(R,u)$, the system of equations characterizing the entropies
$$
q_\rho = u \, \eta_\rho + \rho^{-1} \, \eta_u,
\qquad q_u = \rho \, \eta_{\rho} + u \, \eta_u
$$
reads, by setting $Q := q - u \, \eta$,
\be
Q_R = \eta_u, \qquad Q_u = \eta_R - \eta.
\label{4.10}
\ee
It is clear that the entropy flux can be deduced from the entropy by integration
in $R$ and $u$. We focus attention on the region $0 \leq \rho \leq 1$, that is, $R \leq 0$.
We will use the notation
$$
a \vee b := \max(a,b).
$$

\begin{theorem} {\rm (Entropy-flux kernel.)}
The entropy flux kernel has the form
$$
\sigma (R,u,s) = u\, \chi(R,u-s) + h(R,u-s),
$$
where the function $h$ admits the following representation formulas:
$$
h = - \sgn(u-s) + \frac{\del}{\del u} \int_0^R \chi(r, u-s) \, dr,
$$
or equivalently
\be
h = \frac{\del}{\del s} H(|u-s|,R),
\quad
H = |u-s| + \int_{-(|R| \vee |u-s|)}^{-|u-s|} e^{r/2} \, f(|u-s|^2 - r^2) \, dr,
\label{4.11}
\ee
or still
\be
\begin{split}
h = & \sgn(u-s) \, \bigl( e^{-|u-s|/2} \, \1_{|u-s|<|R|} - 1 \bigr)
   \\
& - 2 \, \int_{-(|R|\vee |u-s|)}^{-|u-s|} (u-s) \, e^{r/2} \, f'(|u-s|^2 - r^2) \, dr.
\end{split}
\label{4.12}
\ee
\end{theorem}

\begin{proof} In view of \eqref{4.10}
we can calculate any value $Q_\ast = Q(R_\ast, u_\ast)$ via an integral, as follows
$$
Q_\ast = \int_{l_\ast} \eta_u \, dR  + (\eta_R -\eta) \, du,\quad
l_\ast = l^-\cup l_0,
$$
where $l^-$ and $l_0$ are the curves in the $(R,u)$--plane given by
\be
\begin{split}
& l^- :  R = 0,   \quad u = \lambda,   \quad -\infty < \lambda < u_\ast,
\\
& l_0 :  R = \lambda R_\ast, \quad u = u_\ast, \quad 0 < \lambda < 1.
\end{split}
\nonumber
\ee
It follows from \eqref{4.9} that
$$
Q_\ast = -2 \, \int_{-\infty}^{u_\ast} \psi (u) \, du
          + \int_0^{R_\ast} \eta_u (R,u_\ast) \, dR.
$$
Substituting $l^-$ by
$$
l^+ : u = \lambda, \quad u_\ast < \lambda <\infty, \quad R = 0,
$$
one obtains similarly that
$$
Q_\ast = - 2 \, \int_\infty^{u_\ast}\psi (u)du
+\int_0^{R_\ast}\eta_u (u_\ast,R)dR.
$$
Observe, that
$$
\int_{-\infty}^{u_\ast} \psi (u) \, du + \int_\infty^{u_\ast} \psi (u) \, du
=
\int_\RR \psi(u) \sgn(u_\ast-u) \, du.
$$
Hence,
$$
Q(R,u) = - \int_\RR \psi(s) \sgn(u-s) \, ds
+
\int_\RR \psi (s) \, \frac{\del}{\del u} \int_0^R\chi (r, u-s) \, dr ds,
$$

Next, we have
$$
\int_0^{R}\chi (r,u-s) \, dr
=
-\int_{-|R|}^0 e^{r/2}f(|u-s|^2 -r^2)\1_{r<-|u-s|}
\1_{r>-|R|}dr = -H_1,
$$
where $H_1$ is the last integral in \eqref{4.11} and, therefore,
the first formula is established.

To calculate
$$
\frac{\del}{\del u}H_1  =
\frac{\del}{\del u}
\int_{-(|R|\vee |u-s|)}^{-|u-s|} e^{r/2} \, f(|u-s|^2 - r^2) \, dr,
$$
we observe that
$$
\frac{\del}{\del u}(|R|\vee |u-s|) = \1_{|u-s|>|R|} \, \sgn(u-s).
$$
Hence, we have
\be
\begin{split}
\frac{\del}{\del u} H_1
= &
2 \, \int_{-(|R|\vee |u-s|)}^{-|u-s|} (u-s) \, e^{r/2} \, f'(|u-s|^2 - r^2) \, dr
-f(0)e^{-|u-s|/2} \sgn(u-s)
\\
& + \1_{|u-s|\geq |R|}
e^{-(|R|\vee |u-s|)/2} \, f(|u-s|^2 - (|R|\vee |u-s|)^2) \, \sgn(u-s).
\end{split}
\nonumber
\ee
The last term coincides  with
$$
\1_{|u-s|\geq |R|} \, e^{-|u-s|/2} \, \sgn(u-s).
$$
Thus, the representation formula \eqref{4.12} is proved and the proof of Theorem 4.5
is completed.
\end{proof}

\


\subsection{Singularities of entropy and entropy-flux kernels}

From the above results we see that the functions $\chi$ and $h$ are
 continuous everywhere
except along the boundary of their support, that is, the lines $u = s \pm |R|$.
The most singular parts (measures and BV part) of the first order derivatives of the
functions $\chi$ and $h$ with respect to the variable $s$ are now computed.

\begin{theorem} {\rm (Singularities of the entropy kernels.)}
The derivatives $\chi_s$ and $h_s$ in $\DD'(\RR)$ are as follows:
\be
\chi_s = e^{R/2} \, \Big( \delta_{s=u-|R|} - \delta_{s=u+|R|} \Big)
   + G^\chi(R,u-s) \, \1_{|u-s|<|R|},
\label{4.13}
\ee
\be
h_s = e^{R/2} \, \Big( \delta_{s=u-|R|} + \delta_{s=u+|R|} \Big)
       + G^h(R,u-s) \, \1_{|u-s|<|R|},
\label{4.14}
\ee
where, for all $|v|\leq |R|$,
$$
G^\chi(R,v) := 2 e^{R/2} \, v \, f'(v^2 -R^2)
$$
and
$$
G^h(R,v) := e^{-|v|/2} \, (1/2 - 2|v|)
     +  2 \, \int_{-(|R|\vee |v|)}^{-|v|} \Big( e^{r/2} \, f'(v^2 -r^2) + 2 e^{r/2} \, v^2 \, f''(v^2 -r^2)
\Big) \, dr.
$$
\end{theorem}

It will be convenient to extend the functions $G^\chi$ and $G^h$
by continuity outside the region $|v|<|R|$ by setting
$$
G^\chi(R,v)
=
\begin{cases}
\hskip.3cm  2 |R| \, f'(0) \, e^{R/2},   & v \geq |R|,
\\
 -2 |R| \, f'(0) \, e^{R/2}, & v \leq -|R|,
\end{cases}
$$
and
$$
G^h(R,v) = e^{R/2} \, (2R + 1/2),  \qquad|v| \geq |R|.
$$

\begin{proof}
Given a test function $\varphi=\varphi(s)$, we can write
\be
\begin{split}
\int_\RR \chi \, \varphi'(s) \, ds
& =
e^{R/2} \, \int_{u-|R|}^{u+|R|}\varphi'(s) \, f(|u-s|^2 -R^2) \, ds
\\
& = e^{R/2} \bigl( \varphi (u+|R|) - \varphi (u-|R|) \bigr)
  + e^{R/2} \int\limits_{u-|R|}^{u+|R|}2\varphi  f'(|u-s|^2 -R^2) (u-s) ds,
\end{split}
\nonumber
\ee
which yields the first formula \eqref{4.13}.

Next, it follows from \eqref{4.11} that
$$
h_s
=\sgn(u-s) \, e^{-|u-s|/2} \, \Big( \delta_{s=u-|R|} - \delta_{s=u+|R|}
  + {1  \over 2} \sgn(u-s)  \1_{|u-s|<|R|}  \Big)
$$
$$
 - 2 \int_{-(|R|\vee |u-s|)}^{-|u-s|} \frac{\del}{\del s}\Big( (u-s) \, e^{r/2}
    \, f'(|u-s|^2 -r^2) \Big) \, dr + 2 \, f'(0) \, e^{-|u-s|/2} \, (u-s) \, \sgn(s-u)
$$
$$
  -2 \, e^{-(|R|\vee |u-s|)/2} \,  \1_{|u-s|\geq |R|} \,
      f'\bigl(|u-s|^2 - (|R|\vee |u-s|)^2\bigr) \, (u-s) \, \sgn(u-s).
$$

The last term above coincides with
$$
-2 \,  e^{- |u-s|/2} \, \1_{|u-s|\geq |R|} \,f'(0)|u-s|
$$
and, therefore, the second formula \eqref{4.14} is also established.
\end{proof}


\section{Reduction of the support of the Young measure}
\setcounter{equation}{0}

\subsection{Tartar's commutation relations}

We now turn to investigating Tartar's commutation relation for Young measures,
following the approach in Chen and LeFloch \cite{CL1,CL2}.
In the previous section we constructed the class of weak entropies $\eta$ and entropy fluxes $q$
in terms of  the variables $\rho$ and $u$. We can also express $\eta$ and $q$ as
functions of the Riemann invariants $W$ and $Z$, via the following change of variables
$$
\bar{\eta }(W,Z) := \eta (u,\rho), \quad \bar{q}(W,Z)=q(u,\rho),
$$
$$
W := \rho e^u, \quad Z := \rho e^{-u}.
$$
To simplify notations, it is convenient to adopt the following convention.
In the rest of this section we will write
$\langle F\rangle = \int F(u,\rho) \, d\nu$ instead of $\int\bar{F}(W,Z) \, d\nu$,
by assuming that $\rho,u$ are the functions of the variables $W,Z$ given by
$$
\rho =(WZ)^{1/2}, \quad
u = \frac{1}{2} \, \ln\frac{W}{Z}.
$$

We will prove:

\begin{theorem} {\rm (Reduction of the support of the Young measure.)}
Let $\nu=\nu(W,Z)$ be a probability measure with support included
in the region
$$
\{(W,Z):0\leq W\leq W_2, \quad0\leq Z\leq Z_2\}
$$
and such that
\be
\la  \eta_1\, q_2 - \eta_2 \,q_1 \ra \,
= \la  \eta_1\ra\, \la  q_2\ra - \la  \eta_2\ra\, \la q_1\ra
\label{5.1}
\ee
(where $\la F\ra :=  \la \nu, F\ra$)
for any two weak entropy pairs $(\eta_1, q_1)$ and $(\eta_2, q_2)$ of the
Euler equations $(1.1)$-$(1.2)$.
Then, the support of $\nu$ in the $(W,Z)$-plane is either a single point or
a subset of the vacuum line $\big\{\rho = 0 \big\} = \{WZ=0\}$.
\end{theorem}

The proof of Theorem 5.1 will be based on {\sl cancellation properties}
associated with the entropy and entropy-flux pairs of systems of conservation laws.
The key idea (going back to DiPerna \cite{DiPerna2})
is that, nearby the diagonal $\bigl\{s_2=s_3\bigr\}$,
the function
$$
E(\rho,v; s_2,s_3): = \chi(\rho, v-s_2) \, \sigma(\rho, v,s_3)
- \chi(\rho, v-s_2) \, \sigma(\rho, v,s_3)
$$
is much more regular than the kernels $\chi$ and $\sigma$ themselves.

The principal scheme can be explained as follows. Given functions
$\psi_i \in$
$\DD(\RR)$, $(i=1,2,3)$,
we define the entropy pairs
$$
\eta_i (u,R)=\int\chi (u-s_i, R)\psi_i (s_i)ds_i, \quad
q_i (u,R)=\int\sigma (u,R,s_i )\psi_i (s_i)ds_i,
$$
and deduce from Tartar's relations \eqref{5.1} the following remarkable identity
(see Chen and LeFloch \cite{CL1}, as well as the earlier work \cite{LPS})
$$
\la \eta_1 q_2 - q_1 \eta_2 \ra \, \la \eta_3 \ra
+ \la q_1 \eta_3 - \eta_1 q_3 \ra \, \la \eta_2 \ra
+ \la q_3 \eta_2 - \eta_3 q_2 \ra \, \la \eta_1 \ra = 0.
$$
Next, substituting $\psi_i (s_i)$ with $-\psi_i'(s_i) = -P_i \psi_i (s_i)$
and denoting $F_i = F(u,R,s_i)$,
we arrive, after cancellation of the arbitrary functions $\psi_i (s_i)$, at the equality
\be
\begin{split}
\la\chi_1 P_2 h_2  -h_1 P_2\chi_2  \ra
\la P_3\chi_3
\ra +
\la h_1 P_3\chi_3 - & \chi_1 P_3 h_3 \ra\la P_2\chi_2 \ra
\\
& = - \la P_3 h_3 P_2 \chi_2 - P_3\chi_3 P_2 h_2 \ra
\la\chi_1\ra,
\end{split}
\label{5.2}
\ee
which is valid in $\DD'(\RR)^3$.
In view of the expression of the distributional derivative of $\sigma$ and $h$ (Theorem~4.6),
each term in \eqref{5.2} can be calculated explicitly.
Denoting
$$
w=u-|R|,\quad z=u+|R|,
$$
we find
\be
\chi_1 P_2 h_2  -h_1P_2 \chi_2  =
e^{R/2}(h_1 -\chi_1)\delta_{s_2 =w}-
e^{R/2}(h_1 +\chi_1)\delta_{s_2 =z}+
(h_1 G^\chi_2-\chi_1 G^h_2)\1_{|u-s_2|<|R|}
\label{5.3}
\ee
and, similarly,
\be
\chi_1 P_3 h_3  -h_1P_3\chi_3  =
e^{R/2}(h_1 -\chi_1)\delta_{s_3 =w}-
e^{R/2}(h_1 +\chi_1)\delta_{s_3 =z}+
(h_1 G^\chi_3-\chi_1 G^h_3)\1_{|u-s_3|<|R|}.
\label{5.4}
\ee
Moreover, we have
\be
\begin{split}
 P_3\chi_3 P_2 h_2  -P_3 h_3  P_2\chi_2
=
&
2e^{R}\Big(\delta_{s_2 =z}\,\delta_{s_3 =w} -\delta_{s_2 =w} \,
\delta_{s_3 =z}\Big)
\\
& + e^{R/2}
\Big(\delta_{s_2 =w}(G_3^\chi-G_3
^h)+\delta_{s_2 =z}(G_3^\chi+G_3
^h )\Big)
\1_{|u-s_3|<|R|}
\\
& + e^{R/2}
\Big(\delta_{s_3 =w}(G_2^h-G_2
^\chi)-\delta_{s_3 =z}(G_2^h+G_2
^\chi)\Big)
\1_{|u-s_2|<|R|}
\\
& + (G_3^\chi G_2^h -G_2^\chi G_3^h)
\1_{|u-s_2|<|R|}\1_{|u-s_3|<|R|}.
\end{split}
\nonumber
\ee

In view of the formulas \eqref{5.3} and \eqref{5.4}
the right-hand side of \eqref{5.2} contains products of functions
with bounded variation (involving $\sigma$ and $h$) and Dirac masses plus smoother terms.
Such products were earlier discussed  by Dal~Maso, LeFloch,
and Murat \cite{DLM}. On the other hand, the right-hand side
of \eqref{5.2} is {\sl more singular} and
involves products of measures, product of BV functions by measures,
and smoother contributions; see \eqref{5.2}.
Our calculation below will show that the left-hand side
of \eqref{5.2} tends to zero in the sense of distributions
if $s_2 \to s_1$
and $s_3 \to s_1$,
while the right-hand side
tends to a (possibly) non-trivial limit.

We test the equality \eqref{5.2} with the function
\be
\psi (s) \, \varphi_2^\eps(s-s_2) \, \varphi_3^\eps(s-s_3)
:=
\psi (s) \, \frac{1}{\eps^2} \, \varphi_2 (\frac{s-s_2}\eps) \, \varphi_3 (\frac{s-s_3}\eps)
\label{5.5}
\ee
of the variables $s = s_1,s_2,s_3$, where $\psi \in \DD(\RR)$ and
$\varphi_j :\RR\to\RR$ is a molifier such that
$$
\varphi_j(s_j) \geq 0, \qquad \int_\RR \varphi_j(s_j) \, ds_j =1, \qquad
\supp \, \varphi_j (s_j) \subset (-1,1).
$$


\subsection{Nonconservative products}

To provide testing of equality  \eqref{5.2} by the function \eqref{5.5}, we will need the following technical
observations.

\begin{lemma}
Let $\psi, F : \RR \to \RR$ and $f : [a',b'] \to \RR$
be continuous functions. Then, for every interval $[a,b] \subset \RR$,
the integral
$$
I^\eps(a,b,a', b') := \int_{a'}^{b'}
\psi (s_1) \, f(s_1) \, \varphi_2^\eps(s_1 -a)
\int_a^b F(s_3) \, \varphi_3^\eps(s_1 -s_3) \, ds_3 ds_1
$$
has the following limit when $\eps \to 0$
$$
\psi(a) \, F(a) \, \Big( A_{2,3}^- \, f(a) \, \1_{a'<a <b'}
+ B_{2,3}^- \, f(a'+) \, \1_{a=a'}
+ C_{2,3}^- \, f(b'-) \, \1_{a=b'} \Big),
$$
where $A_{2,3}^- := B_{2,3}^- + C_{2,3}^-$ and
the coefficients $B^-$ and $C^-$ depend only on the mollifying functions:
$$
B_{2,3}^- := \int_0^\infty \int_{-\infty}^{y_1} \varphi_2 (y_1) \,
\varphi_3 (y_3 ) \, dy_3 dy_1,
\quad
C_{2,3}^- := \int_{-\infty}^0 \int_{-\infty}^{y_1}\varphi_2 (y_1) \,
\varphi_3 (y_3) \, dy_3 dy_1.
$$
\end{lemma}

{\sl Formally} the integral  $I^\eps$ has the form
$$
I(a,b,a', b') := \int_{s_1=a'}^{b'}
\psi (s_1) \, f(s_1) \, \delta_{s_1=a} \int_{s_3=a}^b F(s_3) \, \delta_{s_1 =s_3}.
$$
Lemma 5.2 shows that this term can not be defined in a classical manner and that,
by regularization of the Dirac masses, different limits may be obtained, depending
the choice of the mollifying functions.

Similarly we have

\begin{lemma}
Let $\psi, F : \RR \to \RR$ and $f : [a',b'] \to \RR$
be continuous functions. Then, for every interval $[a,b] \subset \RR$,
the integral
$$
J^\eps(a,b, a', b') := \int_{a'}^{b'}
\psi (s_1) \, f(s_1) \, \varphi_2^\eps(s_1 -b)
\int_{a}^{b} F(s_3) \, \varphi_3^\eps(s_1 - s_3) \, ds_3 ds_1
$$
has the following limit when $\eps \to 0$:
$$
\psi (b) \, F(b) \, \Big(
A_{2,3}^+ \, f(b) \, \1_{a'<b <b'}
+
B_{2,3}^+ \, f(a'+) \, \1_{a=a'}
+ C_{2,3}^+ \, f(b'-) \, \1_{a=b'} \Big),
$$
where $A_{2,3}^+ := B_{2,3}^+ + C_{2,3}^+$ and
$$
B_{2,3}^+ := \int_0^\infty \int_{y_1}^\infty \varphi_2 (y_1) \, \varphi_3 (y_3 ) \, dy_3 dy_1,
C_{2,3}^+
:=
\int_{-\infty}^0 \int_{y_1}^\infty \varphi_2 (y_1) \, \varphi_3 (y_3) \, dy_3 dy_1.
$$
\end{lemma}

\

Along the same lines we have also:

\begin{lemma}
Let $\psi,F:\RR \to \RR$ be continuous functions and
let the function $f :\RR \to \RR$ be continuous
everywhere except possibly at two points $a$ and $b$ with $a<b$. Then,
for every real $\alpha$ the integral
$$
K^\eps(a,b,\alpha) := \int_\RR \psi (s_1) \, f(s_1) \, \varphi_3^\eps(s_1 -\alpha)
\int_a^b  F(s_2) \, \varphi_2^\eps(s_1 -s_2) \, ds_2 ds_1
$$
has the following limit when $\eps \to 0$:
$$
\psi (\alpha) \, F(\alpha) \, \Big( f(\alpha) \, \1_{a<\alpha <b}
+ \big( C_{2,3}^- \, f(a-) + B_{2,3}^- f(a+) \big) \,
\1_{\alpha =a}
+ \big( C_{2,3}^+ \, f(b-) + B_{2,3}^+ \, f(b+) \big) \,
\1_{\alpha =b} \Big).
$$
\end{lemma}

We only give the proof of this last statement. Lemma 5.2 and 5.3 can be checked
similarly.

\

\begin{proof} Making first the change of variables $s_2 = s_1 - \eps y_2$
and then $s_1 = \eps y_1 +\alpha$, one can write
$$
K^\eps= -\int_\RR \psi (\eps y_1 +\alpha) \, f(\eps y_1 +\alpha) \,
\varphi_3 (y_1) \int_{y_1 + (\alpha -a)/ \eps}^{y_1 + (\alpha-b)/\eps} \, F(\eps (y_1 -y_2)+\alpha) \,
\varphi_2 (y_2) \, dy_2 dy_1.
$$
Clearly, we have
$
K^\eps \to 0 \quad \mbox{ when } \alpha <a \mbox{ or } \alpha >b.
$

Now, if $\alpha =a$ we can write
$$
K^\eps = \sum_1^3 K_i^\eps
=
-(\int_{-\infty}^{a}+\int_{a}^{b}+\int_{b}^\infty)
\psi (s_1)f(s_1)\varphi_3^\eps(s_1 -a)
\int_{\frac{s_1 -a}\eps}^{\frac{s_1 -b}\eps}
F(s_1 -\eps y_2)\varphi_2 (y_2) \, dy_2
ds_1.
$$
Consider the first term $K_1^\eps$:
$$
K_1^\eps\stackrel{s_1 =\eps y_1 +a}{=-}
\int_{-\infty}^0\psi (\eps y_1 +a)f(\eps y_1
+a)\varphi_3 (y_1)\int_{y_1}^{y_1 -\frac{b-a}\eps}
F(\eps (y_1 -y_2)+a)\varphi_2 (y_2) \, dy_2 dy_1,
$$
which satisfies
$$
K_1^\eps\to \psi (a) \, F(a)
f(a-)\int_{-\infty}^0\varphi_3 (y_1) \int_{-\infty}^{y_1}
\varphi_2 (y_2) \, dy_2 dy_1.
$$
Similarly, one see that
$$
K_2^\eps\to \psi (a)F(a)
f(a+)\int_0^\infty\varphi_3 (y_1)\int_{-\infty}^{y_1}
\varphi_2 (y_2 ) dy_2 dy_1, \quad
K_3^\eps\to 0.
$$
The other values of $\alpha$ can be studied by the same arguments
and this completes the proof of Lemma 5.4.
\end{proof}


\subsection{Proof of Theorem 5.1}

\

\noindent{\bf Step 1. }
First, we consider the right-hand side of \eqref{5.2}.
Let us denote
$$
d\nu :=d\nu (W,Z),\quad d\nu ':=d\nu (W', Z' ),\quad
w'=u' -|R'|,\quad z'=u'+|R'|
.
$$
Applying the distribution
$\la P_3 h_3 P_2 \chi_2 - P_3\chi_3 P_2 h_2 \ra
\la\chi_1\ra$
to the test function
\eqref{5.5}, we write the integral
\be
\int_{\RR^3}  \la P_3 h_3(s_3)  P_2 \chi_2(s_2) - P_3\chi_3(s_3) P_2 h_2(s_2) \ra \,
\la \chi_1(s) \ra \, \psi(s) \, \varphi_2^\eps(s - s_2) \,
\varphi_3^\eps(s - s_3) \, ds ds_2 ds_3
\label{5.6}
\ee
as the sum
$
 \sum_1^4 I_i^\eps,
$
in which, in view of Theorem 4.6, we can distinguish between
products of Dirac measures
$$
I_1^\eps :=
\int \psi\la\chi_1\ra
\Big\la
2e^{R}\varphi_2^\eps(s -z))\varphi_3^\eps(s -w))-
2e^{R}\varphi_2^\eps(s -w))\varphi_3^\eps(s -z))
\Big\ra ds,
$$
products of Dirac measure by functions with bounded variation
\be
\begin{split}
I_2^\eps := &
\int \psi\la\chi_1\ra
\Big\la e^{R/2}\varphi_2^\eps(s -w))
\int (G_3^\chi-G_3^h)\varphi_3^\eps(s -
s_3)\1_{|u-s_3|<|R|}ds_3
\Big\ra ds
\\
& -
\int \psi\la\chi_1\ra \Big\la
e^{R/2}\varphi_3^\eps(s -w))
\int (G_2^\chi-G_2^h)\varphi_2^\eps(s -s_2)\1_{|u-s_2|<|R|}ds_2
\Big\ra ds
\\
& = : I_{2,1}^\eps - I_{2,2}^\eps,
\end{split}
\nonumber
\ee
\be
\begin{split}
I_3^\eps :=&
\int \psi\la\chi_1\ra
\Big\la
e^{R/2}\varphi_2^\eps(s -z))
\int (G_3^\chi+G_3^h)\varphi_3^\eps(s -s_3)\1_{|u-s_3|<|R|} \, ds_3
\Big\ra \, ds
\\
&
- \int \psi\la\chi_1\ra
\Big\la
e^{R/2}\varphi_3^\eps(s -z))
\int (G_2^\chi+G_2^h)\varphi_2^\eps(s -s_2)\1_{|u-s_2|<|R|} \, ds_2
\Big\ra ds,
\end{split}
\nonumber
\ee
and a smoother remainder
$$
I_4^\eps=
\int\limits_{\RR^3} \psi\la\chi_1\ra
\Big\la
(G_3^\chi G_2^h-G_2^\chi G_3^h)\1_{|u-s_3|<|R|} \, \1_{|u-s_2|<|R|}
\Big\ra \varphi_2^\eps(s -s_2)\varphi
_3^\eps(s -s_3) \, ds ds_2 ds_3.
$$

By change of variable we see that the integral
\be
\begin{split}
I_1^\eps =
{2 \over \eps}  \int\limits_{W,Z} \int\limits_{W',Z'}
\int e^ R & \, \psi(\eps y + z) \, \chi (R', u' - (\eps y + z))
\\
& \Big( \varphi_2(y) \, \varphi_3(y + 2|R|/\eps)
        - \varphi_2 (y + 2|R|/\eps) \, \varphi_3(y)\Big) \, dy d\nu d\nu ' .
\end{split}
\nonumber
\ee
tends to zero : $I_1^\eps\to 0$.
The same is true for the smoothest term $I_4^\eps$, in view of the identity
\be
\begin{split}
I_4^\eps
= \int \int \int \Big(&
\int_{w}^{z} G^\chi(R,u,s_3) \, \varphi_3^\eps(s -s_3) \, ds_3
\int_{w}^{z} G^h(R,u,s_2) \, \varphi_2^\eps(s - s_2) \, ds_2
\\
& -\int_{w}^{z} G^\chi(R,u,s_2) \, \varphi_2^\eps(s - s_2) \, ds_2
\int_{w}^{z} G^h(R,u,s_3) \, \varphi_3^\eps(s - s_3) \, ds_3
\Big)
\\
& \hskip5.5cm \psi \chi (R', u' - s) \, ds d\nu d\nu ',
\end{split}
\nonumber
\ee
which clearly tends to
\be
\begin{split}
\int \int \int \psi \, \chi (R', u' - s)
\, \1_{|u-s|<|R|} \, \Big( & G^\chi(R,u,s) \, G^h(R,u,s)
\\
& - G^\chi(R,u,s) \, G^h(R,u,s) \Big) \, d s d\nu d\nu ' = 0.
\end{split}
\nonumber
\ee

\

We denote
$$
Q^\pm := G^{\chi}\pm G^h, \quad
F_i = F(R,u,s_i),\quad F_i ' = F(R',u',s_i).
$$

Let us now consider the term $I_2^\eps = I_{2,1}^\eps - I_{2,2}^\eps$
in \eqref{5.6}.
We have
$$
I_{2,1}^\eps = \int \int\int \psi e^{R/2} \chi_1'
\, \varphi_2^\eps(s -w) \int_{w}^{z}
Q_3^- \, \varphi_3^\eps(s -s_3) \, ds_3 ds d\nu d\nu'.
$$
Therefore, in view of Lemma 5.2, we obtain that $I_{2,1}^\eps$ tends to
$$
\iint  e^{R/2}  \psi (w)  Q^- (w)
  \Big( \chi'(w)  \1_{w'< w <z'}  A_{2,3}^- + \chi'(w'+)
 \1_{w=w'}  B_{2,3}^-+  \chi'(z'-) \1_{w=z'}  C_{2,3}^- \Big)
 d\nu d\nu',
$$
and $I_{2,2}^\eps$ tends to
$$
\iint e^{R/2} \psi (w) \, Q^- (w)
\Big( \chi'(w) \1_{w'< w <z'} A_{3,2}^-
+ \chi'(w'+) \1_{w=w'}B_{3,2}^- + \chi'(z'-)\1_{w=z'} C_{3,2}^-\Big)
d\nu d\nu',
$$
as $\eps\to 0$.
We conclude that the limit of $I_2^\eps$ is equal to
$$
\iint e^{R/2} \, \psi (w) \, Q^-(w)
\, \Big( \chi'(w) \, \1_{w'<w<z'} \, A^-
 + \chi'(w'+) \, \1_{w=w'} \, B^- + \chi'(z'-) \, \1_{w=z'} \, C^- \Big)
\, d\nu d\nu',
$$
where
$$
A^-=A_{2,3}^--A_{3,2}^-,\quad
B^-=B_{2,3}^--B_{3,2}^-,\quad C^-=C_{2,3}^--C_{3,2}^-.
$$

By Lemma 5.3 we can determine similarly that $\lim_{\eps \to 0} I_3^\eps $ is equal to
$$
\iint  e^{R/2} \, \psi (z) \, Q^+ (z)
\, \Big( \chi'(z) \, \1_{w'<z<z'}A^+
   + \chi'(w'+) \, \1_{z=w'} \, B^+ + \chi'(z'-) \, \1_{z=z'} \, C^+ \Big)
d\nu d\nu',
$$
where
$$
A^+ := A_{2,3}^+ - A_{3,2}^+,
\quad
B^+ := B_{2,3}^+ - B_{3,2}^+,
\quad C^+ := C_{2,3}^+ - C_{3,2}^+.
$$
In conclusion we have identified the limit of the term \eqref{5.10}, it is equal to
$$ 
\iint  e^{R/2} \psi (w) Q^- (w)
\Big( \chi'(w) \1_{w'<w<z'}  A^- + \chi'(w'+)  \1_{w=w'}  B^-
+ \chi'(z'-)  \1_{w=z'}  C^- \Big)  d\nu d\nu'
$$
$$
+ \iint e^{R/2}  \psi(z)  Q^+(z)
\Big( \chi'(z)  \1_{w'<z<z'}  A^+
 + \chi'(w'+)  \1_{z=w'}B^+ + \chi'(z'-)  \1_{z=z'}  C^+ \Big)
 d\nu d\nu'. 
$$


\

\noindent{\bf Step 2. }
We now proceed by studying the two terms in the left-hand side of \eqref{5.2}.
We apply the distribution
$\la\chi_1(s_1) P_2 h_2(s_2) - h_1(s_1) P_2 \chi_2(s_2) \ra
   \, \la P_3\chi_3(s_3) \ra$
to the test function
\eqref{5.5}.
We write the integral
$$
\int_{\RR^{3}} \la\chi_1(s_1) P_2 h_2(s_2) - h_1(s_1) P_2 \chi_2(s_2) \ra
   \, \la P_3\chi_3(s_3) \ra \,
   \psi(s_1) \, \varphi_2^\eps(s_1 -s_2) \, \varphi_3^\eps(s_1 -s_3) \, ds_1 ds_2ds_3
$$
as the sum
$
\sum_1^3 J_i^\eps,
$
where
$$
J_1^\eps := \int_{\RR^3} \la e^{R/2} \, (h_1 - \chi_1)
\, \delta_{s_2 =w} \ra \,
\la P_3 \chi_3 \ra \, \psi (s) \, \varphi_2^\eps(s -s_2) \,
\varphi_3^\eps(s -s_3) \, ds ds_2 ds_3
$$
$$
J_2^\eps := - \int_{\RR^3} \la e^{R/2}(h_1 +\chi_1)
\delta_{s_2 =z} \ra \, \la P_3\chi_3\ra
\, \psi (s) \, \varphi_2^\eps(s -s_2) \, \varphi_3^\eps(s -s_3) \, ds ds_2 ds_3,
$$
and
$$
J_3^\eps :=
\int_{\RR^3}\la (h_1 G^\chi_2-\chi_1 G^h_2)\1_{|u-s_2|<|R|}
\ra
\la P_3\chi_3\ra\psi (s)\varphi_2^\eps(s -s_2)\varphi
_3^\eps(s -s_3) \, ds ds_2 ds_3 .
$$

The application of the distribution
$
\la \chi_1 P_3 h_3^\eps - h_1 P_3 \chi_3^\eps \ra
\la P_2\chi_2^\eps \ra
$
to the test function
\eqref{5.5} can be represented similarly. The integral
$$
\int_{\RR^3}
\la \chi_1(s) P_3 h_3(s_3) - h_1(s) P_3 \chi_3(s_3) \ra
\,
\la P_2\chi_2(s_2) \ra \, \psi(s) \, \varphi_2^\eps(s -s_2)
\, \varphi_3^\eps(s -s_3) \, ds ds_2 ds_3
$$
is equal to
$\sum_1^3 K_i^\eps$,
where
$$
K_1^\eps := \int_{\RR^3} \la e^{R/2} \, (h_1 -\chi_1) \, \delta_{s_3 =w} \, \ra
\la P_2\chi_2\ra \, \psi (s) \, \varphi_2^\eps(s -s_2) \, \varphi_3^\eps(s -s_3) \, ds ds_2 ds_3,
$$
$$
K_2^\eps := - \int_{\RR^3}\la e^{R/2} \, (h_1 + \chi_1) \,
\delta_{s_3 =z} \ra \, \la P_2 \chi_2\ra \,
\psi (s)\varphi_2^\eps(s -s_2) \, \varphi_3^\eps(s -s_3) \, ds ds_2 ds_3,
$$
and
$$
K_3^\eps := \int_{\RR^3} \la (h_1 \, G^\chi_3 - \chi_1 G^h_3) \, \1_{|u-s_3|<|R|}\ra \,
   \la P_2\chi_2\ra \, \psi (s) \, \varphi_2^\eps(s -s_2) \,
   \varphi_3^\eps(s -s_3) \, dsds_2 ds_3 .
$$

Some further decomposition of these integral terms will be necessary:
$$
J_i^\eps = \sum_1^3 J_{i,j}^\eps,
\quad
K_i^\eps = \sum_1^3 K_{i,j}^\eps,
$$
where
$$
J_{1,1}^\eps := \int_{\RR^3}\la e^{R/2}(h_1 -\chi_1)
\delta_{s_2 =w} \ra \la e^{R/2} \delta_{s_3 =w}
\ra\quad\psi (s_1)\varphi_2^\eps(s -s_2)\varphi
_3^\eps(s -s_3)  \, ds ds_2 ds_3,
$$
$$
J_{1,2}^\eps := - \int_{\RR^3}
\la e^{R/2}(h_1 -\chi_1)
\delta_{s_2 =w}
\ra
\la e^{R/2}
\delta_{s_3 =z}
\ra\quad\psi (s_1)\varphi_2^\eps(s -s_2)\varphi
_3^\eps(s -s_3) \, ds ds_2 ds_3,
$$
$$
J_{1,3}^\eps := \int_{\RR^3}
\la e^{R/2}(h_1 -\chi_1)
\delta_{s_2 =w}
\ra
\la
G_3^\chi\1_{|u-s_3 |<|R|}
\ra\quad\psi (s)\varphi_2^\eps(s -s_2)\varphi
_3^\eps(s -s_3) \, ds ds_2 ds_3,
$$
$$
J_{2,1}^\eps :=- \int_{\RR^3}
\la e^{R/2}(h_1 +\chi_1)
\delta_{s_2 =z}
\ra
\la e^{R/2}
\delta_{s_3 =w}
\ra\quad\psi (s)\varphi_2^\eps(s -s_2)\varphi
_3^\eps(s -s_3) \, ds ds_2 ds_3,
$$
$$
J_{2,2}^\eps := \int_{\RR^3}
\la e^{R/2}(h_1 +\chi_1)
\delta_{s_2 =z}
\ra
\la e^{R/2}
\delta_{s_3 =z}
\ra\quad\psi (s)\varphi_2^\eps(s -s_2)\varphi
_3^\eps(s -s_3) \, ds ds_2 ds_3,
$$
$$
J_{2,3}^\eps =- \int_{\RR^3}
\la e^{R/2}(h_1 +\chi_1)
\delta_{s_2 =z}
\ra
\la
G_3^\chi\1_{|u-s_3 |<|R|}
\ra\quad\psi (s)\varphi_2^\eps(s -s_2)\varphi
_3^\eps(s -s_3) \, ds ds_2 ds_3,
$$
$$
J_{3,1}^\eps := \int_{\RR^3}
\la
(h_1 G^\chi_2-\chi_1 G^h_2)\1_{|u-s_2|<|R|}
\ra
\la
e^{R/2}
\delta_{s_3 =w}
\ra\quad\psi (s)\varphi_2^\eps(s -s_2)\varphi
_3^\eps(s -s_3) \, ds ds_2 ds_3,
$$
$$
J_{3,2}^\eps :=- \int_{\RR^3}
\la
(h_1 G^\chi_2-\chi_1 G^h_2)\1_{|u-s_2|<|R|}
\ra
\la
e^{R/2}
\delta_{s_3 =z}
\ra\quad\psi (s)\varphi_2^\eps(s -s_2)\varphi
_3^\eps(s -s_3) \, ds ds_2 ds_3,
$$
$$
J_{3,3}^\eps := \int_{\RR^3}
\la
(h_1 G^\chi_2-\chi_1 G^h_2)\1_{|u-s_2|<|R|}
\ra
\la
G_3^\chi\1_{|u-s_3 |<|R|}
\ra\quad\psi (s)\varphi_2^\eps(s -s_2)\varphi
_3^\eps(s -s_3) \, ds ds_2 ds_3 .
$$
The terms $K_{1,1}^\eps$, $K_{1,2}^\eps$, etc. are defined in a completely
analogous fashion.

We can put $J_{1,1}^\eps$ in the form
\be
\begin{split}
J_{1,1}^\eps
= & {1 \over \eps^2} \iiint  \psi \, e^{(R + R')/2} \, (h_1 -\chi_1)
\, \varphi_2 (\frac{s -w)}\eps)
\varphi_3 (\frac{s -w')}\eps) \, d\nu d\nu' ds
\\
= & {1 \over \eps}
\iiint  \psi (s) e^{(R+ R')/2} \, (h -\chi )(s)
\, \varphi_2 (y_1)|_{s =\eps y_1 +w}
\varphi_3 (y_1 +\frac{w -w')}\eps)
d\nu d\nu' dy_1.
\end{split}
\nonumber
\ee
A similar representation formula is valid for $K_{1,1}^\eps$. In consequence we find
\be
\begin{split}
J_{1,1}^\eps - K_{1,1}^\eps
= & \frac{1}{\eps }
\iiint  \psi (s) \, e^{(R + R')/2} \, (h -\chi )(s) |_{s =\eps y_1+w}
\\
& \Big(\varphi_2 (y_1) \, \varphi_3( y_1 + \frac{w- w')}
\eps)
- \varphi_3 (y_1) \, \varphi_2(y_1 +\frac{w-w')} \eps) \Big) d\nu d\nu' dy_1.
\end{split}
\nonumber
\ee
Clearly, we have
$$
J_{1,1}^\eps - K_{1,1}^\eps \to 0.
$$
The terms $J_{k,l}^\eps$ and  $K_{k,l}^\eps$ contain the product
of measures or the product of BV-functions and can be treated in the
same manner. In turn, one obtains
$$
J_{1,2}^\eps - K_{1,2}^\eps \to 0,
\,
J_{2,1}^\eps - K_{2,1}^\eps \to 0,
\,
J_{2,2}^\eps - K_{2,2}^\eps \to 0,
\, J_{3,3}^\eps - K_{3,3}^\eps \to 0.
$$

Let us consider the terms $J_{k,l}^\eps$ and $K_{k,l}^\eps$, containing
the product of a measure and a BV-function. By Lemma 5.4, the term
$$
J_{1,3}^\eps =
\iiint  \psi (s)e^{R/2} \, (h_1 -\chi_1) \, \varphi_2^\eps(s -w)
\int_{w'}^{z'} {G_3^\chi}' \, \varphi_3^\eps(s-s_3) \, ds_3 ds d\nu
d\nu'
$$
converges toward
$$
\iint  e^{R/2} \, \psi (w) \, (h-\chi)(w) \, {G^\chi}'(w)
\Big(  \1_{w'<w<z'} + \1_{w=w'} \, (C_{2,3}^- + B_{2,3}^-)
+ \1_{w=z'} \, (C_{2,3}^+ + B_{2,3}^+ )\Big)  d\nu
d\nu',
$$
hence,
\be
\begin{split}
& \lim_{\eps \to 0} (J_{1,3}^\eps - K_{1,3}^\eps)
\\
& =  \iint  e^{R/2} \, \psi (w) \, (h-\chi)(w) \, {G^\chi}'(w)
\Big(  \1_{w=w'} \, (C^-+B^-)+ \1_{w=z'}(C^+ +B^+ )
\Big)d\nu d\nu'.
\end{split}
\nonumber
\ee

By the same argument we find that the term
$$
J_{3,1}^\eps
=
\iiint    \psi (s) \, e^{R'/2} \varphi_3^\eps(s -w')
\, \int_{w}^{z} (h_1 G_2^\chi - \chi_1 G_2^h ) \, \varphi_2^\eps(s
-s_2) \, ds_2 dsd\nu d\nu'
$$
tends toward
$$
\iint  e^{R'/2} \, \psi (w') \, G^\chi(w') \, L \, d\nu d\nu'
-
\iint  e^{R'/2} \, \psi (w') \, G^h(w') \, S \, d\nu d\nu',
$$
where
$$
L :=  \1_{w<w'<z} \, h(w') +
\1_{w'=w} \, \Big( h(w-) \, C_{3,2}^- + h(w+) \, B_{3,2}^- \Big)
+
\1_{w'=z} \, \Big( h(z-) \, C_{3,2}^+ + h(z+) \, B_{3,2}^+ \Big)
$$
and
$$
S :=  \1_{w<w'<z} \, \chi(w') + \1_{w'=w} \, \Big( \chi (w-) \, C_{3,2}^-
 + \chi (w+) \, B_{3,2}^-\Big) + \1_{w'=z} \, \Big( \chi (z-)C_{3,2}^+ +\chi (z+)B_{3,2}^+ \Big).
$$
Hence,
$$
J_{3,1}^\eps-K_{3,1}^\eps\to
-\iint  e^{R'/2} \psi (w') \, (\1_{w'=w} M_w + \1_{w'=z} \, M_z) \, d\nu
d\nu',
$$
where
$$
M_w := G^\chi(w)\Big( h(w-)C^-+h(w+)B^-\Big) -G^h (w)
\Big( \chi (w-)C^- +\chi (w+)B^- \Big),
$$
$$
M_z =
G^\chi(z)\Big( h(z-)C^+ +h(z+)B^+ \Big) -G^h (z)
\Big( \chi (z-)C^+ +\chi (z+)B^+ \Big).
$$
One can see
that the integrals, containing the functions
$\1_{w'=w}$ and $\1_{w'=z}$ cancel each other.

In a similar way, we can treat the other terms and arrive at the final
equality
$$
\iint e^{R/2} \, \psi(w) \, Q^-(w) \, \chi'(w) \, A^- \,
\1_{w'<w<z'} +  e^{R/2} \, \psi(z) \, Q^-(z) \, \chi'(z) \, A^+ \,
\1_{w'<z<z'} \, d\nu d\nu' = 0,
$$
resulting from
\eqref{5.2}.

\


\noindent{\bf Step 3. }
Observing that $A^+ =-A^-$ and
$$
Q^-(w)=-Q^+ (z)=
e^{R/2}\Big( -\frac{f(0)}2+2|R|+2|R|f'(0)\Big) = : D(R),
$$
we can write
\be
\iint Y \, e^{R/2} \, D(R)
\Big( \psi(w) \, \chi'(w) \, \1_{w'<w<z'}
+
\psi(z) \, \chi'(z) \, \1_{w'<z<z'}
\Big)  \, d\nu d\nu'=0,
\label{5.7}
\ee
where
$$
Y=\int\limits_{-\infty}^{+\infty}\int\limits_{-\infty}^{s_2}
\varphi_2 (s_2)\varphi_3 (s_3)-\varphi_2 (s_3)\varphi_3 (s_3)ds_2 ds_3 .
$$
As observed in \cite{CL2}, the functions $\varphi_2 $ and $\varphi_3$
can easily be chosen in such a way that $Y \ne 0$.

In accordance with our notation, the equality \eqref{5.7} means precisely
\be
\begin{split}
\iint & D(\frac{1}{2}\ln (WZ)) \, (W\, Z)^{1/4} \,
( W' \, Z' )^{1/4}
\, \Big(  f\bigl(( - \ln{W' \over W} \, \ln (Z' W) \bigr) \, \psi (\ln W) \,
\1_{W' < W < 1/Z'}
\\
& + f\bigl(- \ln (W' Z) \, \ln {Z' \over Z}\bigr) \, \psi (-\ln Z) \,
\1_{ W' < 1/Z < 1/Z'}
\Big)  \, d\nu  d\nu '  = 0.
\end{split}
\nonumber
\ee
The conditions \eqref{3.20} guarantee that $|D(R)|\geq e^{R/2}/2$.
Since $f(-x^2) \geq 1$ and the function $\psi$ is arbitrary, we conclude from the last equality that
$$
\iint_{W,Z} (WZ)^{1/2} \iint_{\bigl\{W' < W\bigr\} \cap \bigl\{Z'<1/W\bigr\}}
(W' Z')^{1/4} \, d\nu (W',Z')d\nu (W,Z)=0
$$
and
$$
\iint_{W,Z} (WZ)^{1/2} \iint_{ \bigl\{W' < 1/Z \bigr\} \cap \bigl\{ Z'<Z \bigr\}}
(W'Z')^{1/4} \, d\nu (W',Z') \, d\nu (W,Z) = 0.
$$

We arrive at the following important claim: whenever
$$
W^\ast Z^\ast \ne 0, \quad
(W^\ast,Z^\ast) \in \supp \nu,
$$
we have
$$
\iint_{W,Z} (WZ)^{1/2} d\nu (W,Z) \neq 0
$$
and therefore
\be
\begin{split}
& \iint_{ \bigl\{ W' < W^\ast\bigr\} \cap \bigl\{Z'<1/W^\ast \bigr\}}
(W'Z')^{1/4} \, d\nu (W',Z') = 0,
\\
& \iint_{\bigl\{W'<1/Z^\ast \bigr\} \cap \bigl\{ Z'<Z^\ast \bigr\}}
(W'Z')^{1/4} \, d\nu (W',Z') = 0.
\end{split}
\label{5.8}
\ee

We will now conclude from \eqref{5.8} that the Young measure is a Dirac mass
or a measure concentrated at the vacuum.
At this stage, it is useful to draw a picture on the $W,Z-$plane, with the $W-$axis
being horizontal.
One should draw two hyperbolas $WZ=1$ and $WZ=\rho_2^2$, keeping in mind that
$\mbox{supp}\, \nu$ lies below the hyperbola
$WZ=\rho_2^2$, where the constant $\rho_2 < 1$ is defined in section 3.
The hyperbola $WZ=1$ helps to picture the set
$$
M^\ast := \Big(
\bigl\{0< W' < W^\ast \bigr\} \cap \bigl\{0< Z' < 1/W^\ast\bigr\}
\Big)
\cup
\Big( \bigl\{ 0<W' < 1/Z^\ast \bigr\} \cap \bigl\{0< Z' < Z^\ast \bigr\}
\Big),
$$
a union of two rectangulars. The relations \eqref{5.8}
imply that  $M^\ast$ does not
intersect the support of $\nu$ :
\be
W^\ast Z^\ast \ne 0\quad\mbox { and } \quad
(W^\ast,Z^\ast) \in \supp \nu\quad\Longrightarrow
 M^\ast \cap \supp \nu = \emptyset.
\label{5.9}
\ee
By construction, the hyperbola
$WZ=1$ does not intersect  $\mbox{supp}\, \nu $. The inclusion
$\mbox{supp}\, \nu \subset \{\rho =0\}$ holds
if no hyperbola
$WZ=\delta$, $0<\delta  <1$, intersects  $\mbox{supp}\, \nu $.
If $\mbox{supp}\, \nu $ contains a point $(W,Z)$ such that
$\rho (W,Z)\ne 0$, there is a number $0< \delta <1$ such that the hyperbola
$WZ=\delta$, intersects $\mbox{supp}\, \nu $.
Let $0<\delta_0 <1$ be the largest number such that
the hyperbola $WZ=\delta_0$ intersects  $\mbox{supp}\, \nu $.
By \eqref{5.9}, the intersection
$$
\mbox{supp}\, \nu \cap \{WZ=\delta_0\}
$$
may contain only one point $(W^{\ast},Z^{\ast})$ and
$$
\mbox{supp}\, \nu \cap \{ 0<WZ <\delta_0\}=\emptyset .
$$
Thus
\be
\nu =\alpha \delta_{\ast}+\mu,
\label{5.10}
\ee
with $\mbox{supp}\, \mu\subset \{\rho =0\}$.
Throughout the paper we use only weak entropies.
Hence, putting \eqref{5.10} into the Tartar's commutation relation
\eqref{5.1}, we obtain that any two entropy pairs satisfy the equality
\be
\alpha (q_2 \eta_1 -q_1 \eta_2)=\alpha^2 (q_2 \eta_1 -q_1 \eta_2)
\label{5.11}
\ee
at the point $(W^{\ast},Z^{\ast})$. Let us choose the following
entropy pair (as in \cite{Shelukhin2})
$$
\eta_i =\rho^{B_i}e^{A_i u},\quad
q_i =-\frac{ A_i}{B_i -1}\rho^{B_i -1}e^{A_i u},\quad A_i =\sqrt {B_i (B_i -1)},
\quad B_1 \ne B_2.
$$
Now, the equality\eqref{5.11} is rewritten as
$$
\alpha (1-\alpha)\rho_{\ast}^{B_1 +B_2 -1}e^{(A_1 +A_2)u_{\ast}}
\Big( \sqrt {\frac{B_1}{B_1 -1}}-\sqrt {\frac{B_2}{B_2 -1}}\Big)=0.
$$
Hence, $\alpha =0$ or $\alpha =1$. This completes the proof of Theorem 5.1.


\section{Convergence and compactness of solutions}
\setcounter{equation}{0}

Due to the decomposition \eqref{5.10} of the Young measures, the convergence formulas
\eqref{3.9} imply that
$$
W^\eps \rightharpoonup W,\quad Z^\eps \rightharpoonup Z,\quad
F(W^\eps,Z^\eps) \rightharpoonup F(W,Z)\quad \mbox{ weakly $\star$ in } L^\infty(\Pi),
$$
for any function $F(\alpha, \beta)$, $F\in C(K)$, such that $F=0$ at the vacuum set $\alpha\beta=0$.
(See formula \eqref{3.19} for the definition of the compact set $K$.)
Hence, for almost all $(x,t)\in \Pi$
$$
\rho^\eps :=
(W^\eps Z^\eps)^{1/2}\to \rho =(WZ)^{1/2} =: f_1 (W,Z),
$$
$$
m^\eps := (W^\eps \, Z^\eps)^{1/2}\ln
(\frac{W^\eps}{Z^\eps})^{1/2}\to
m =
(WZ)^{1/2}\ln (\frac{W}{Z})^{1/2}
=: f_2 (W,Z),
$$
$$
\frac{(m^\eps)^2}{\rho^\eps} :=
(W^\eps \, Z^\eps)^{1/2}(\ln(\frac{W^\eps}{Z^\eps})^{1/2})^2\to
\frac{m^2}{\rho}
= f_3 (W,Z) := (WZ)^{1/2}(\ln (\frac{W}{Z})^{1/2})^2.
$$
Moreover,
\be
F(m^\eps,\rho^\eps)\to
F(m,\rho)\quad \mbox{for almost all}\quad (x,t)\in \Pi
\label{6.1}
\ee
for any function $F(m,\rho)$ such that
\be
\tilde{F}(\alpha, \beta) :=  F(f_2(\alpha,\beta),f_1
(\alpha,\beta))\in C(K), \quad \tilde{F}|_{\alpha\beta =0}=0.
\label{6.2}
\ee
Indeed, one can derive the convergence
\eqref{6.1} from the following fact:
$$
v^\eps\to v \quad\mbox{ and } \quad (v^\eps)^2 \to v^2 \quad
\mbox{ weakly in } L_{loc}^2(\Pi) \Longrightarrow v^\eps\to v\quad
\mbox{ strongly in } L_{loc}^2(\Pi).
$$

Let us show that $(m,\rho)$ is an entropy solution of problem
\eqref{2.1a}. To this end we let $\eps$ and $\eps_1$ go to zero in
\eqref{3.12}. (More exactly we should do it in the similar equality relevant to the
auxiliary approximation.)
If  functions $\eta (m,\rho)$, $q(m,\rho)$ obey the restrictions
\eqref{6.2}, one obtains
$$
\int(\eta (m^\eps,\rho^\eps)-
\eta (m_0^\eps,\rho_0^\eps))\varphi_t +
q(m^\eps,\rho^\eps)\varphi_x \, dxdt\to
\int (\eta (m,\rho )
- \eta (m_0, \rho_0))\varphi_t + q(m,\rho )\varphi_x \, dxdt,
$$
$$
\eps\int\eta (m^\eps,\rho^\eps)\varphi_{xx} \, dxdt\to 0
$$
for any
$\varphi\in\DD(\RR^2)$.

From now on we assume that $\eps_1 =\eps^r$, $r>1$.
If a function $\eta (m,\rho)$ meets the conditions of Theorem 2.1, the derivatives
$\eta_m (m,\rho)$ and $m_{\rho} (m,\rho)$ are continuous on any closed set
$$
\{0\leq\rho\leq\rho_1, \quad
|m|\leq c_1 \rho (1+|\ln\rho|)\},\quad \rho_1 >0,\quad c_1 >0 .
$$
Hence, by estimate \eqref{3.18},
$$
|\eps_1 u_x (q_m +\eta_{\rho})|=
|2\eps_1 u_x (\frac{m}{\rho}\eta_m +\eta_{\rho})|
\leq c\eps_1 (|uu_x| +|u_x|)\leq
\eps^{1/2}\rho^{1/2}|u_x|
(\eps^{\frac{r-1}{2}}+|u|\rho^{\gamma}\eps
^{\delta}),
$$
where
$2\gamma <\frac{r-1}{r}$ $,2\delta =r(1-2\gamma)-1$.
Besides,
$
\eps_1 \rho^{-1}|\eta_m \rho_x|\leq
$
$
c\eps
^{1/2}\rho^{-1/2}|\rho_x|\eps^{\frac{r-1}{r}}.
$
Now, it follows from Lemma 3.3 and estimates \eqref{3.18} that
$$
\eps_1 u_x (q_m +\eta_{\rho})
-2\eps_1
\eta_m \rho^{-1}\rho_x \to 0
\quad\mbox{in}\quad L_{loc}^2 (\Pi).
$$
Taking into account the convexity of the function $\eta (m,\rho)$, we send $\eps$ to zero in
\eqref{3.12} to deduce that the pair $(m,\rho)$ is an entropy
solution of \eqref{2.1a}-\eqref{2.1b}. The proof of Theorems 2.1 to 2.3 is completed.

We conclude by giving a proof of Theorem 2.4.
Let $(m_n, \rho_n)$ be a sequence of bounded in $L^\infty(\Pi)$ entropy solutions of the problem
\eqref{2.1a} obeying the restriction of Theorem 2.4. We introduce the sequences
$$
W_n := \rho_n e^{m_n /\rho_n}, \quad Z_n = \rho_n \, e^{-m_n /\rho_n}.
$$
Clearly, we have
$$
W_n \rightharpoonup W,\quad
Z_n \rightharpoonup Z\quad \mbox{ weakly $\star$ in } L_{loc}^\infty(\Pi),
$$
and there exist Young measures $\nu_{x,t}$ such that, for all $F(\alpha,\beta) \in C_{loc}(\RR^2)$,
$$
F(W_n (x,t),Z_n (x,t)) \rightharpoonup \la\nu_{x,t},F\ra.
$$
Given two entropy pairs $(\eta_i (m,\rho),q_i (m,\rho))$ from Theorem 2.1,
the sequences of measures
$$
\theta_i^n :=
\del_t \eta_i (m_n, \rho_n)+\del_x q_i (m_n, \rho_n)=
\del_t \tilde{\eta}(W_n, Z_n)+\del_x \tilde{q}(W_n, Z_n),
$$
satisfy the conditions of Murat's lemma and are compact in $W_{loc}^{-1,2}(\Pi)$.
(We remind the notation
$\tilde{q}(W,Z):=$
$q(f_2 (W,Z),f_1 (W,Z))$.)
By the $div-curl$ lemma, the Tartar commutation relations \eqref{3.14}
is valid for the Young measures $\nu_{x,t}$. Then we argue like in the proof of Theorem 2.1 to arrive at
the decomposition \eqref{5.10} for $\nu_{x,t}$.
Hence,
$$
F(W_n (x,t),Z_n (x,t))\to F(W(x,t),Z(x,t))\quad\mbox{ almost everywhere in}\quad \Pi
$$
for any
$ F(\alpha,\beta)\in C_{loc}(\RR^2)$ such that $F(\alpha,\beta)=0$ if $\alpha\beta =0$.
Denoting
$$
\rho =f_1 (W,Z),\quad m=f_2 (W,Z),\quad
\rho_n =f_1 (W_n, Z_n),\quad
m_n =f_2 (W_n, Z_n),
$$
we pass to the limit, as $n\to\infty$, in the inequality
$$
\iint  \eta(m_n, \rho_n) \, \del_t \varphi + q(m_n, \rho_n) \, \del_x \varphi \quad dxdt
+ \int \eta(m_0, \rho_0) \, \varphi(x,0) \, dx \geq 0
$$
and check that $(m,\rho)$ is an entropy solution of the problem
\eqref{2.1a}-\eqref{2.1b}. The proof of Theorem 2.4 is completed.


\section*{Acknowledgments}
The authors were supported by a grant from INTAS (01-868).
The support and hospitality of the Isaac Newton Institute
for Mathematical Sciences, University of Cambridge, where part of this research was performed
during the Semester Program ``Nonlinear Hyperbolic Waves in Phase Dynamics and Astrophysics''
(January to July 2003) is also gratefully acknowledged.
P.G.L. was also supported by the Centre National de la Recherche
Scientifique (CNRS).


\

\

\newcommand{\auth}{\textsc}
\newcommand{\tit}{\textrm}
\newcommand{\jou}{\textit}


\begin{thebibliography}{10}

\bibitem{Ball} \auth{Ball J.M.,}
A version of the fundamental theorem for Young measures.
In: ``PDE's and Continuum Models of Phase Transitions'', Lecture Notes
in Physics, Vol. 344, Rascle M., Serre D., and Slemrod M. (eds),
Springer Verlag,, pp.~241--259.

\bibitem{BB} \auth{Bianchini S. and Bressan A.,} 
Vanishing viscosity solutions of nonlinear hyperbolic systems,
{\em Ann. of Math.} (2004). 

\bibitem{Chen} \auth{Chen G.-Q.,}  
Convergence of the Lax-Friedrichs scheme for isentropic 
gas dynamics (III), Acta Math. Sci. 8 (1988), 243--276. 

\bibitem{CL1} \auth{Chen G.-Q. and LeFloch P.G.,}
Compressible Euler equations with general pressure law,
Arch. Rational Mech Anal. 153 (2000), 221--259.

\bibitem {CL2} \auth{Chen G.-Q. and LeFloch P.G.,}
Existence theory for the isentropic Euler equations,
Arch. Rational Mech. Anal. 166 (2003), 81--98.

\bibitem{Dafermos} \auth{Dafermos C.M.,}
{\sl Hyperbolic conservation laws in continuum physics,}
Grundlehren Math. Wissenschaften Series 325, Springer Verlag, 2000.

\bibitem{DLM} \auth{Dal Maso G., LeFloch P.G., and Murat F.,}
Definition and weak stability of nonconservative products,
J. Math. Pures Appl. 74 (1995), 483--548.

\bibitem{DCL}  \auth{Ding X., Chen G.-Q., and Luo P.,}
Convergence of the Lax-Friedrichs scheme for the isentropic gas dynamics,
Acta Math. Sci. 5 (1985), 483--540.

\bibitem{DiPerna1} \auth{DiPerna R.J.,}
Convergence of approximate solutions to conservation laws,
Arch. Rational Mech. Anal. 82 (1983), 27--70.

\bibitem{DiPerna2} \auth{DiPerna R.J.,}
Convergence of the viscosity method for isentropic gas dynamics,
Commun. Math. Phys. 91 (1983), 1--30.

\bibitem{FS} \auth{Frid H. and Shelukhin V.V.,}
A quasilinear parabolic system for three-phase capillary flow in porous media,
SIAM J. Math. Anal. 35 (2003), 1029--1041.

\bibitem{HW} \auth{Huang F.-M. and Wang Z.,}
Convergence of viscosity solutions for isothermal
gas dynamics, SIAM J. Math. Anal. 34 (2003), 595--610.

\bibitem{LSU} \auth{Lady\v{z}enskaja O.A., Solonnikov V.A., and Ural'ceva N.N.,}
{\sl Linear and quasi-linear equations of parabolic type,}
A.M.S., Providence, 1968.

\bibitem{Lax} \auth{Lax P.D.,}
{\sl Hyperbolic systems of conservation laws and the mathematical theory of shock waves,}
SIAM Regional Conf. Lecture 11, Philadelphia, 1973.

\bibitem{LeFloch1} \auth{LeFloch P.G.,}
Existence of entropy solutions for the compressible Euler equations,
International Series Numer. Math. Vol.~130, Birk\"auser Verlag B\"asel, Switzerland,
1999, pp.~599--607.

\bibitem{LeFloch2} \auth{LeFloch P.G.,}
{\sl Hyperbolic systems of conservation laws:
The theory of classical and nonclassical shock waves},
Lectures in Mathematics, ETH Z\"urich, Birk\"auser, 2002.

\bibitem{LeFloch3} \auth{LeFloch P.G.,}  
to appear. 

\bibitem{LPT} \auth{Lions P.L., Perthame B., and Tadmor E.,}
Kinetic formulation for the isentropic gas dynamics and p-system,
Commun. Math. Phys. 163 (1994), 415--431.

\bibitem{LPS} \auth{Lions P.L., Perthame B., and Souganidis P.E.,}
Existence and stability of entropy solutions for the hyperbolic
systems of isentropic gas dynamics in Eulerian and Lagrangian coordinates,
Comm. Pure Appl. Math. 49 (1996), 599--638.

\bibitem{Morawetz} \auth{Morawetz C.J.,}
An alternative proof of DiPerna's theorem,
Comm. Pure Appl. Math. 44 (1991), 1081--1090.

\bibitem{Lu} \auth{Lu Y.-G.,}
Convergence of the viscosity method for nonstrictly hyperbolic conservation laws,
Comm. Math. Phys. 150 (1992), 59--64.

\bibitem{Murat1} \auth{Murat F.,}
Compacit\'e par compensation,
Ann. Scuola Norm. Sup. Pisa Sci. Fis. Mat. 5 (1978), 489--507.

\bibitem{Murat2} \auth{Murat F.,}
L'injection du c\^one positif de  $H^{-1}$
dans $W^{-1, q}$ est compacte pour tout $q<2$,
J. Math. Pures Appl. 60 (1981), 309--322.

\bibitem{Nishida} \auth{Nishida T.,}
Global solutions for the initial-boundary value problem
of a quasilinear hyperbolic systems,
Proc. Japan Acad. 44 (1968), 642--646.

\bibitem{Olver} \auth{Olver F.W.J.,}
{\sl Asymptotics and special functions,}
Academic Press, 1974.

\bibitem{Perthame} \auth{Perthame B.,}
{\sl Kinetic formulation of conservation laws,}
Lecture Series in Math. and Appl., Oxford Univ. Press, 2002.
             
\bibitem{PT} \auth{Perthame B. and Tzavaras A.,}
Kinetic formulation for systems of two conservation laws  and elastodynamics,
Arch. Rational Mech. Anal. 155 (2000), 1--48.

\bibitem{Ovsyannikov} \auth{Ovsyannikov L.V.,}
{\sl Group analysis of differential equations,} Academic Press, New York, 1982.

\bibitem{Serre} \auth{Serre D.,}
La compacit\'e par compensation pour les syst\`emes hyperboliques non-lin\'eaires 
de deux \'equations \`a une dimension d'espace,
J. Math. Pures Appl. 65 (1986), 423--468.

\bibitem{Shelukhin1} \auth{Shelukhin V.V.,}
Existence theorem in the variational problem for
compressible inviscid fluids, Manuscripta Math. 61 (1988), 495--509.

\bibitem{Shelukhin2} \auth{Shelukhin V.V.,}
Compactness of bounded quasientropy solutions to the system of equations of an isothermal gas.
Siberian Math. J. 44 (2003), 366--377.

\bibitem{Tartar1} \auth{Tartar L.,}
Compensated compactness and applications to partial differential equations,
in ``Nonlinear analysis and mechanics: Heriot-Watt Symposium'', Vol.~IV,
Res. Notes in Math., 1979, Vol.~39, Pitman, Boston, Mass.-London, pp.~136--212.

\bibitem{Tartar2} \auth{Tartar L.,}
The compensated compactness method applied to systems of conservation laws,
in ``Systems of Nonlinear Partial Differential Equations'', J.M. Ball ed.,
NATO ASI Series, C. Reidel publishing Col., 1983, pp.~263--285.

\end{thebibliography}
\end{document}